\documentclass{amsproc}
\usepackage{amssymb}
\usepackage{color}
\usepackage{units}   

\numberwithin{equation}{section}

\newtheorem{thm}{Theorem}[section]
\newtheorem{cor}[thm]{Corollary}
\newtheorem{lem}[thm]{Lemma}
\newtheorem{pro}[thm]{Proposition}
\newcounter{step}

\theoremstyle{remark}
\newtheorem{rem}[thm]{Remark}
\newtheorem{exa}[thm]{Example}

\theoremstyle{definition}

\DeclareMathOperator{\supp}{supp}

\newcommand*{\card}{\mathrm{card\,}}
\newcommand*{\cbb}{\mathbb C}
\newcommand*{\D}{\mathrm d\,}

\newcommand*{\E}{\mathrm{e}}

\newcommand*{\Ge}{\geqslant}

\newcommand*{\Le}{\leqslant}

\newcommand*{\nbb}{\mathbb N}
\newcommand*{\nfr}{\mathfrak n}

\newcommand*{\rbb}{{\mathbb R}}

\newcommand*{\sqrtk}{\sqrt[\leftroot{1}\uproot{4}\varkappa]}
\newcommand*{\sqrtx}[2]{\sqrt[\leftroot{1}\uproot{#1}#2]}
\newcommand*{\sfr}{\mathfrak S}

\newcommand*{\ssqrt}[1]{\sfr_{#1}^{\nicefrac{1}{2}}}
\newcommand*{\ssqrtt}[2]{\sfr_{#1,#2}^{\nicefrac{1}{2}}}

\newcommand*{\zbb}{\mathbb Z}

\begin{document}
   \title[On the $\varkappa\,$th root of a Stieltjes moment
sequence] {On the $\varkappa\,$th root of a Stieltjes
moment sequence}
   \author[J.\ Stochel]{Jan Stochel}
   \address{Instytut Matematyki, Uniwersytet Jagiello\'nski,
ul.\ {\L}ojasiewicza 6, PL-30348 Kra\-k\'ow,
Poland}
   \email{Jan.Stochel@im.uj.edu.pl}
   \author[J.\ B.\ Stochel]{Jerzy Bart{\l}omiej Stochel}
\address{Faculty of Applied Mathematics,
AGH University of Science and
Technology, Al. Mickiewicza 30, PL-30059 Krak\'ow,
Poland}
   \email{stochel@agh.edu.pl}
  \thanks{This work was  supported by the
MNiSzW grant No.\ NN201 546438 (2010-2013) as well as
by the MNiSzW grant No.\ 10.42004 (the second
author).}

   \keywords{Stieltjes moment sequence,
$\varkappa\,$th root of a Stieltjes moment sequence,
representing measure, support of a measure, hole of a
support}
   \begin{abstract}
Stieltjes moment sequences $\{a_n\}_{n=0}^\infty$ whose
$\varkappa\,$th roots $\{\sqrtk{a_n}\}_{n=0}^\infty$
are Stieltjes moment sequences are studied ($\varkappa$
is a fixed integer greater than or equal to $2$). A
formula connecting the closed supports of representing
measures of $\{a_n\}_{n=0}^\infty$ and
$\{\sqrtk{a_n}\}_{n=0}^\infty$ is established. The
relationship between the holes of the supports of these
measures is investigated. The set of all pairs $(M,N)$
of positive integers for which there exists a Stieltjes
moment sequence whose square root is a Stieltjes moment
sequence and both of them have representing measures
supported on subsets of $(0,\infty)$ of cardinality $M$
and $N$, respectively, is described.
   \end{abstract}
   \maketitle
   \section{Introduction}
In \cite{Hor0,Hor1,Hor2} Horn laid the foundations of
the theory of infinitely divisible matrices, kernels
and positive definite sequences. Among other things,
he described Stieltjes moment sequences
$\{a_n\}_{n=0}^\infty$ whose all powers
$\{a_n^\alpha\}_{n=0}^\infty$ with positive real
exponents $\alpha$ are Stieltjes moment sequences
(cf.\ \cite[Theorem 2.9]{Hor2}). In a recent paper
\cite{jsjs} we asked the question:\ for what positive
integers $M$ is it true that the square root of a
Stieltjes moment sequence $\big\{\alpha_1
\vartheta_1^n + \ldots + \alpha_M
\vartheta_M^n\big\}_{n=0}^\infty$ is not a Stieltjes
moment sequence for all real numbers $\alpha_1,
\ldots, \alpha_M > 0$ and $0 <\vartheta_1 < \ldots <
\vartheta_M$? We will show that the answer to this
question is in the affirmative for\footnote{\;The case
of $M=2$, which was answered in \cite[Lemma
3.3]{jsjs}, immediately implies that Rhoades' square
root problem (cf.\ \cite[p.\ 296]{Rh}) has a negative
solution (see \cite[Theorem 7]{Ben} for a solution of
this problem based on a particular choice of a
Hausdorff moment sequence whose representing measure
is supported on a two point set).} $M \in \{2,4\}$,
and in the negative for $M \notin \{2,4\}$ (cf.\
Corollary \ref{discr}).

The present work is motivated by the aforementioned
papers (including \cite{Rh,Ben}). We will consider
Stieltjes moment sequences $\{a_n\}_{n=0}^\infty$ whose
$\varkappa\,$th roots $\{\sqrtk{a_n}\}_{n=0}^\infty$
are Stieltjes moment sequences, where $\varkappa$ is a
fixed integer greater than or equal to $2$; by the
Schur product theorem (see \cite[p.\ 14]{sch} or
\cite[Theorem 7.5.3]{Hor-Joh}) this is equivalent to
considering the $\varkappa\,$th powers of Stieltjes
moment sequences. Under the assumption that
$\{a_n\}_{n=0}^\infty$ is determinate, we give a
formula for the (closed) support of a representing
measure $\mu$ of $\{a_n\}_{n=0}^\infty$ written in
terms of the support of a representing measure $\nu$ of
$\{\sqrtk{a_n}\}_{n=0}^\infty$ (see Theorem
\ref{suppmu}). In Section \ref{sec4} and \ref{secnew}
we provide some solutions to the following problem:
given a hole of the support of the measure $\mu$,
determine the circumstances under which the support of
$\nu$ has a hole and then localize it (see Theorems
\ref{lem1+} and \ref{17.02.11}); the converse problem
is studied as well (see Theorem \ref{lem1}). Some
solutions of this problem are written in terms of the
parameters $\iota_{\mathrm{s}}$ and
$\iota_{\mathrm{s}}^*$ that describe, in a sense, the
geometry of the hole of the support of $\mu$ (see
Theorem \ref{17.02.11}). Using the results of Sections
\ref{sec2}-\ref{secnew}, we construct a variety of
examples illustrating the concepts of the paper (see
Section \ref{przyk}). In particular, an example of a
Stieltjes moment sequence which $\varkappa\,$th root is
a Stieltjes moment sequence for $\varkappa=2,4$, but
not for $\varkappa=3$ is furnished (see Example
\ref{prznotsq}). We conclude this paper with a
description of the set of all pairs $(M,N)$ of positive
integers for which there exists a Stieltjes moment
sequence whose square root is a Stieltjes moment
sequence and both of them have representing measures
supported on subsets of $(0,\infty)$ of cardinality $M$
and $N$, respectively (see Theorem \ref{tw2}). The
reader is encouraged to refer to
\cite{akh,b-ch-r,sh-tam,sim} for the mathematical
details of the theory of classical moment problems.
   \section{\label{sec2}Preliminaries}
From now on, the fields of real and complex numbers
are denoted by $\rbb$ and $\cbb$, respectively, and
$\zbb$ stands for the set of all integers. Set
$\rbb_+=\{x \in \rbb\colon x \Ge 0\}$, $\zbb_+ = \{n
\in \zbb\colon n \Ge 0\}$ and $\nbb=\{n \in \zbb\colon
n\Ge 1\}$. Given $x \in \rbb$, we define $\lfloor x
\rfloor = \max\{n \in \zbb \colon n\Le x\}$ and
$\lceil x \rceil= \min\{n \in \zbb\colon x \Le n\}$.
The cardinality of a set $X$ is denoted by $\card
(X)$. $X^{\varkappa}$ stands for the $\varkappa$-fold
Cartesian product of $X$ by itself. We write $\supp
\mu$ for the (closed) support of a regular positive
Borel measure $\mu$ on a Hausdorff topological space.
Given $\theta \in \rbb_+$, we denote by
$\delta_\theta$ the Borel probability measure on
$\rbb_+$ concentrated at $\{\theta\}$.

A sequence $\{a_n\}_{n=0}^\infty \subseteq \rbb_+$ is
said to be a {\em Stieltjes moment} sequence if there
exists a positive Borel measure\footnote{\;Such $\mu$
being finite is automatically regular (see e.g.,
\cite[Theorem 2.18]{rud}).} $\mu$ on $\rbb_+$ such
that $a_n = \int_{\rbb_+} x^n \D \mu(x)$ for all $n
\in \zbb_+$; such $\mu$ is called a {\em representing
measure} of $\{a_n\}_{n=0}^\infty$. If a Stieltjes
moment sequence has a unique representing measure, we
call it {\em determinate}. Recall that (cf.\
\cite{fug})
   \begin{align}\label{compsup}
   \text{\begin{minipage}{70ex} each Stieltjes
moment sequence which has a compactly supported
representing measure is determinate.
   \end{minipage}}
   \end{align}
It is well known that if $(X, \varSigma, \mu)$ is a
measure space, $f$ is a complex $\varSigma$-measurable
function on $X$ and $\int_X |f|^{r} \D \mu < \infty$
for some $r \in (0,\infty)$, then the $\mu$-essential
supremum of $|f|$ is equal to $\lim_{p \to \infty}
\Big(\int_X |f|^p \D \mu\Big)^{1/p}$ (cf.\ \cite[p.\
71]{rud}). This implies that
   \begin{align} \label{esssup}
   \begin{minipage}{28em}
if $\{a_n\}_{n=0}^\infty$ is a Stieltjes moment
sequence with a representing measure $\mu$, then
$ \lim_{n \to \infty} a_n^{1/n} = \sup \supp
\mu$.
   \end{minipage}
   \end{align}

Given an integer $\varkappa \Ge 2$, we define the
continuous mapping $\pi_\varkappa\colon
\rbb_+^\varkappa \to \rbb_+$ by
   \begin{align}     \label{pikap}
\pi_\varkappa (x_1, \ldots, x_\varkappa)=x_1 \cdots
x_\varkappa, \quad x_1, \ldots, x_\varkappa \in
\rbb_+.
   \end{align}
   \begin{lem} \label{pikappa}
Let $F$ be a subset of $\rbb_+$ and $\varkappa$ be an
integer greater than or equal to $2$. If $F$ is
compact, then $\pi_\varkappa(F^\varkappa)$ is compact.
If $0 \notin F$ and $F$ is closed, then $0 \notin
\pi_\varkappa(F^\varkappa)$ and
$\pi_\varkappa(F^\varkappa)$ is closed.
   \end{lem}
   \begin{proof}[Proof of Lemma  \ref{pikappa}]
Clearly, it is enough to consider the case when
$0\notin F$ and $F$ is closed. We claim that if
$\{(x^{(n)}_1, \ldots,
x^{(n)}_\varkappa)\}_{n=1}^\infty \subseteq
F^\varkappa$ and the sequence
$\{\pi_\varkappa(x^{(n)}_1, \ldots,
x^{(n)}_\varkappa)\}_{n=1}^\infty$ is convergent in
$\rbb_+$, then there exists a subsequence of
$\{(x^{(n)}_1, \ldots,
x^{(n)}_\varkappa)\}_{n=1}^\infty$ which is convergent
in~ $F^\varkappa$.

First note that the sequence
$\{x^{(n)}_j\}_{n=1}^\infty$ is bounded for each $j
\in \{1, \ldots, \varkappa\}$. Indeed, otherwise, by
passing to a subsequence if necessary, we may assume
that for some $k \in \{1, \ldots, \varkappa\}$,
$x^{(n)}_k \to \infty$ as $n\to \infty$. Since $F$ is
closed and $0\notin F$, we deduce that $\inf_{n\Ge 1}
x^{(n)}_{j} > 0$ for all $j \in \{1, \ldots,
\varkappa\}$. This implies that the sequence
$\{\pi_\varkappa(x^{(n)}_1, \ldots,
x^{(n)}_\varkappa)\}_{n=1}^\infty$ is unbounded, which
is a contradiction. Now using the compactness
argument, we establish our claim.

By the assertion just proved, we see that $0 \notin
\pi_\varkappa(F^\varkappa)$ and
$\pi_\varkappa(F^\varkappa)$ is closed.
   \end{proof}
Note that Lemma \ref{pikappa} is no longer true if the
assumption $0\notin F$ is dropped. Indeed, if $\varkappa\Ge
2$ and $F=\{0,1,2,\ldots\} \cup \{\frac 12, \frac 13, \frac
14, \ldots\}$, then $F$ is closed, while
$\pi_\varkappa(F^\varkappa)$, being equal to the set of all
nonnegative rational numbers, is not closed.
   \section{The relationship between $\supp\mu$ and $\supp\nu$\label{sec3}}
In the present paper we consider the following
situation:
   \begin{align} \label{zal}
   \begin{minipage}{75ex}
$\{a_n\}_{n=0}^\infty$ is a Stieltjes moment sequence
with a representing measure $\mu$, $\varkappa$ is an
integer greater than or equal to $2$ and
$\{\sqrtk{a_n}\,\}_{n=0}^\infty$ is a Stieltjes moment
sequence with a representing measure $\nu$.
   \end{minipage}
   \end{align}
We begin by proving the basic relation between $\mu$
and $\nu$.
   \begin{lem} \label{nSti}
If {\em \eqref{zal}} holds and $\{a_n\}_{n=0}^\infty$
is determinate, then
   \begin{align}  \label{rhsof}
\mu(E) = \nu^{\otimes
\varkappa}(\pi_\varkappa^{-1}(E)) \quad \text{for all
Borel subsets $E$ of $\rbb_+$,}
   \end{align}
where $\nu^{\otimes \varkappa}$ is the
$\varkappa$-fold product of the measure $\nu$ by
itself and $\pi_\varkappa$ is as in \eqref{pikap}.
   \end{lem}
   \begin{proof}
It follows from Fubini's theorem and \cite[Theorem C,
p.\ 163]{hal2} that
   \begin{multline*}
a_n = \int_{\rbb_+^\varkappa} (x_1 \cdots
x_\varkappa)^n \D \nu^{\otimes \varkappa} (x_1,\ldots,
x_\varkappa) = \int_{\rbb_+} z^n \D \big(\nu^{\otimes
\varkappa} \circ \pi_\varkappa^{-1}\big) (z), \quad
n\in \zbb_+,
   \end{multline*}
where $\nu^{\otimes \varkappa} \circ
\pi_\varkappa^{-1}$ is a positive Borel measure on
$\rbb_+$ given by the right hand side of the equality
in \eqref{rhsof}. Thus, by the determinacy of
$\{a_n\}_{n=0}^\infty$, the measures $\mu$ and
$\nu^{\otimes \varkappa} \circ \pi_\varkappa^{-1}$
coincide.
   \end{proof}
The following lemma which describes the support of the
transport of a measure is surely folklore. For the
reader's convenience, we include its proof.
   \begin{lem} \label{newlem}
Let $X$ and $Y$ be Hausdorff topological spaces,
$\rho$ and $\sigma$ be regular positive Borel measures
on $X$ and $Y$, respectively, and $\phi \colon X \to
Y$ be a continuous mapping such that $\sigma(E) =
\rho(\phi^{-1}(E))$ for all Borel subsets $E$ of $Y$.
Then
   \begin{align*}
\supp \sigma = \overline{\phi(\supp \rho)}.
   \end{align*}
   \end{lem}
   \begin{proof}
We begin by showing that $\phi(\supp \rho) \subseteq
\supp \sigma$. Take $x \in \supp \rho$. Let $V$ be an
open neighbourhood of $\phi(x)$. By the continuity of
$\phi$, the set $\phi^{-1}(V)$ is an open
neighbourhood of $x$, and thus $\sigma(V) =
\rho(\phi^{-1}(V)) > 0$. Since $V$ is an arbitrary
open neighbourhood of $\phi(x)$, we get $\phi(x) \in
\supp \sigma$.

It remains to show that $Y \setminus
\overline{\phi(\supp \rho)} \subseteq Y \setminus
\supp \sigma$. Take $y \in Y \setminus
\overline{\phi(\supp \rho)}$. Then $V:= Y \setminus
\overline{\phi(\supp \rho)}$ is an open neighbourhood
of $y$ such that $V \cap \phi(\supp\rho) =
\varnothing$. This implies that $\phi^{-1}(V) \cap
\supp\rho = \varnothing$. Hence $\sigma(V) =
\rho(\phi^{-1}(V)) = 0$, which yields $y \in Y
\setminus \supp \sigma$. This completes the proof.
   \end{proof}
We are now ready to provide a formula connecting
$\supp\mu$ with $\supp\nu$.
   \begin{thm}\label{suppmu}
Suppose {\em \eqref{zal}} holds and
$\{a_n\}_{n=0}^\infty$ is determinate. Then
   \begin{align}  \label{porfelTer}
\supp \mu = \overline{\pi_\varkappa\big((\supp
\nu)^\varkappa\big)}.
   \end{align}
Moreover, if $\supp \nu$ is compact or $0 \not\in
\supp\nu$, then $\supp \mu = \pi_\varkappa\big((\supp
\nu)^\varkappa\big)$.
   \end{thm}
   \begin{proof}
Applying Lemmata \ref{nSti} and \ref{newlem} and the
well known equality $\supp \nu^{\otimes \varkappa} =
(\supp \nu)^{\varkappa}$, we get \eqref{porfelTer}.
The ``moreover'' part follows from \eqref{porfelTer}
and Lemma \ref{pikappa}.
   \end{proof}
We will show in Example \ref{notclos} that the closure
sign in \eqref{porfelTer} cannot be omitted.
   \begin{cor} \label{suppmu-cor}
If {\em \eqref{zal}} holds and $\{a_n\}_{n=0}^\infty$
is determinate, then
   \begin{enumerate}
   \item[(i)] $\supp \nu \subseteq
\sqrt[\leftroot{1}\uproot{4}\varkappa]{\supp \mu}$,
   \item[(ii)] $\card(\supp \nu) \Le
\card(\supp \mu)$,
   \item[(iii)] $\sup \supp \nu = \sqrtk{\sup \supp
\mu}$.
   \end{enumerate}
   \end{cor}
   \begin{proof}
Conditions (i) and (iii) follow from Theorem
\ref{suppmu} (condition (iii) can also be deduced from
\eqref{esssup}). Condition (ii) is a direct
consequence of (i).
   \end{proof}
As shown in Example \ref{ontheleft}, inclusion (i) in
Corollary \ref{suppmu-cor} may be proper. In turn,
inequality (ii) in Corollary \ref{suppmu-cor} may be
strict (however, in view of Theorem \ref{suppmu}, if
$\card(\supp \mu) = \aleph_0$, then $\card(\supp \nu)
= \card(\supp \mu)$). In fact, it may happen that
$\supp \nu$ is discreet and has only one accumulation
point, while $\supp \mu = \rbb_+$ (cf.\ Example
\ref{notclos}).
   \section{Transforming holes of $\supp\mu$ and $\supp\nu$}
Suppose \eqref{zal} holds. In this and the subsequent
two sections we will study the relationship between
the following two situations:
   \begin{align*}
& \supp \mu \subseteq [0,\vartheta_1] \cup
[\vartheta_2, \vartheta_3] \text{ with } \vartheta_1,
\vartheta_2, \vartheta_3 \in \rbb \text{ such that } 0
\Le \vartheta_1 < \vartheta_2 \Le \vartheta_3,
   \\
& \supp \nu \subseteq [0,\alpha] \cup [\beta, \gamma]
\text{ with } \alpha, \beta, \gamma \in \rbb \text{
such that } 0 \Le \alpha < \beta \Le \gamma.
   \end{align*}
Hereafter we will consider a transformation
$(\vartheta_1, \vartheta_2, \vartheta_3) \to (\alpha,
\beta, \gamma)$ between triplets of real numbers
satisfying the inequalities $0 \Le \vartheta_1 <
\vartheta_2 \Le \vartheta_3$ and $0 \Le \alpha < \beta
\Le \gamma$ which is given by
   \begin{align} \label{theta2.5}
\alpha=\frac{\vartheta_1}{\vartheta_3}
\sqrtk{\vartheta_3}, \quad \beta=\sqrtk{\vartheta_2}
\quad \text{and} \quad \gamma=\sqrtk{\vartheta_3}.
   \end{align}
This transformation is well defined (because
$(\frac{\vartheta_1}{\vartheta_3})^\varkappa <
(\frac{\vartheta_2}{\vartheta_3})^\varkappa \Le
\frac{\vartheta_2}{\vartheta_3}$) and injective, but
not surjective. A triplet $(\alpha, \beta, \gamma)$
with $0 \Le \alpha < \beta \Le \gamma$ is the image of
some $(\vartheta_1, \vartheta_2, \vartheta_3)$ under
this transformation if and only if $\alpha
\gamma^{\varkappa - 1} < \beta^\varkappa$. If this is
the case, then
   \begin{align*}
\vartheta_1 = \alpha \gamma^{\varkappa - 1}, \quad
\vartheta_2=\beta^\varkappa \quad \text{and} \quad
\vartheta_3=\gamma^\varkappa.
   \end{align*}

Given real numbers $\vartheta_1, \vartheta_2,
\vartheta_3$ such that $0 \Le \vartheta_1 <
\vartheta_2 \Le \vartheta_3$, we set
   \begin{align} \label{abplus}
\alpha^\dag=\frac{\vartheta_2}{\vartheta_3}
\sqrtk{\vartheta_3}, \quad \beta^\dag
=\sqrtk{\vartheta_1}.
   \end{align}
The quantities just defined will play an essential
role in Theorem \ref{lem1+} below. Clearly, if
$\alpha$, $\beta$ and $\gamma$ are defined by
\eqref{theta2.5}, then
   \begin{align} \label{duzoal}
\text{$0 \Le \alpha < \beta \Le \gamma$, $\alpha <
\alpha^\dag \Le \beta$, $0 \Le \beta^\dag < \beta$ and
$\alpha \gamma^{\varkappa-1} < \beta^\varkappa$.}
   \end{align}
Note also that
   \begin{align*}
\text{$\alpha^\dag < \beta$ if and only if
$\vartheta_2 < \vartheta_3$.}
   \end{align*}
In general, there is no order relationship between
$\alpha^\dag$ and $\beta^\dag$. Obviously we have
   \begin{align*}
\alpha^\dag < \beta^\dag \iff \vartheta_2^\varkappa <
\vartheta_1 \vartheta_3^{\varkappa-1} \quad \text{and}
\quad \alpha^\dag = \beta^\dag \iff
\vartheta_2^\varkappa = \vartheta_1
\vartheta_3^{\varkappa-1}.
   \end{align*}
The reader should be aware of the fact that the
quantities $\alpha$, $\beta$, $\gamma$, $\alpha^\dag$
and $\beta^\dag$ depend not only on $(\vartheta_1,
\vartheta_2, \vartheta_3)$ but also on $\varkappa$.
Making explicit the dependence on $\varkappa$, we can
formulate the following result.
   \begin{pro} \label{podloga}
If $0 \Le \vartheta_1 < \vartheta_2 \Le \vartheta_3 <
\infty$, then
   \begin{enumerate}
   \item[(i)]
$\lim_{\varkappa \to \infty} \alpha^\dag(\varkappa) =
\vartheta_2/\vartheta_3$,
   \item[(ii)] $\lim_{\varkappa \to \infty} \beta^\dag(\varkappa)=1$
provided that $\vartheta_1 >0$,
   \item[(iii)] if $\alpha^\dag(\varkappa)
< \beta^\dag(\varkappa)$ for some $\varkappa \Ge 2$,
then $\alpha^\dag(\varkappa^\prime) <
\beta^\dag(\varkappa^\prime)$ for all
$\varkappa^\prime \Ge \varkappa$.
   \end{enumerate}
Moreover, if $0 < \vartheta_1 < \vartheta_2 <
\vartheta_3 < \infty$, then
   \begin{enumerate}
   \item[(iv)]
$\alpha^\dag(\iota_{\mathrm{s}}) <
\beta^\dag(\iota_{\mathrm{s}})$,
   \item[(v)] the sequence $\{\beta^\dag(\varkappa) -
\alpha^\dag(\varkappa)\}_{\varkappa =
\iota_{\mathrm{s}}}^\infty$ is strictly increasing
provided that either $\vartheta_1 \Le 1$ or
$\vartheta_1 > 1$ and $\log \vartheta_1 \Le
\frac{\vartheta_2}{\vartheta_3} \log \vartheta_3$,
   \end{enumerate}
where
   \begin{align} \label{iotadef}
\iota_{\mathrm{s}} = \iota_{\mathrm{s}}(\vartheta_1,
\vartheta_2, \vartheta_3) := 1 + \left\lfloor
\frac{\log{(\vartheta_3/\vartheta_1)}}
{\log{(\vartheta_3/\vartheta_2)}}\right\rfloor \Ge 2.
   \end{align}
   \end{pro}
   \begin{proof} The properties (i)-(iv)
are easily seen to be true.

(v) Define the function $f\colon [\iota_{\mathrm{s}},
\infty) \to \rbb_+$ by $f(x) =
\sqrtx{2}{x}{\vartheta_1} -
\frac{\vartheta_2}{\vartheta_3}
\sqrtx{2}{x}{\vartheta_3}$ for $x\in
[\iota_{\mathrm{s}}, \infty)$. It follows from (iv)
that
   \begin{align} \label{strin}
\frac{\vartheta_2}{\vartheta_3} <
\sqrtx{4}{\iota_{\mathrm{s}}}{\frac{\vartheta_1}{\vartheta_3}}
\Le \sqrtx{4}{x}{\frac{\vartheta_1}{\vartheta_3}},
\quad x\in [\iota_{\mathrm{s}}, \infty).
   \end{align}
Note that
   \begin{align} \label{num1}
-x^2 f^\prime(x) = \sqrtx{2}{x}{\vartheta_1} \log
\vartheta_1 - \frac{\vartheta_2}{\vartheta_3}
\sqrtx{2}{x}{\vartheta_3} \log{\vartheta_3}, \quad
x\in [\iota_{\mathrm{s}}, \infty).
   \end{align}
If $\vartheta_1 \Le 1$, then
   \begin{align*} -x^2
f^\prime(x) \overset{\eqref{num1}}<
\Big(\sqrtx{2}{x}{\vartheta_1} -
\frac{\vartheta_2}{\vartheta_3}
\sqrtx{2}{x}{\vartheta_3}\Big) \log{\vartheta_1}
\overset{\eqref{strin}} \Le 0, \quad x\in
[\iota_{\mathrm{s}}, \infty).
   \end{align*}
In turn, if $\vartheta_1 > 1$ and $\log \vartheta_1
\Le \frac{\vartheta_2}{\vartheta_3} \log \vartheta_3$,
then
   \begin{align*} -x^2
f^\prime(x) \overset{\eqref{num1}} \Le
\Big(\sqrtx{2}{x}{\vartheta_1} -
\sqrtx{2}{x}{\vartheta_3}\Big) \log{\vartheta_1} < 0,
\quad x\in [\iota_{\mathrm{s}}, \infty).
   \end{align*}
Hence, in both cases, $f^\prime (x) > 0$ for all $x\in
[\iota_{\mathrm{s}}, \infty)$, which implies that $f$
is strictly increasing. Since $f(\varkappa) =
\beta^\dag(\varkappa) - \alpha^\dag(\varkappa)$ for
$\varkappa=2, 3, \ldots,$ the proof is complete.
   \end{proof}
Regarding Proposition \ref{podloga}\,(v), we note that
the sequence $\{\beta^\dag(\varkappa) -
\alpha^\dag(\varkappa)\}_{\varkappa =
\iota_{\mathrm{s}}}^\infty$ may be strictly decreasing
(e.g., consider $\vartheta_1 = 2$, $\vartheta_2 = 5/2$
and $\vartheta_3 = 2 \cdot 10^4$).
   \section{Relating  holes of $\supp\mu$ and
$\supp\nu$\label{sec4}}
   Our next goal is to analyze the situation when the
support of a representing measure $\nu$ of the
$\varkappa\,$th root of a Stieltjes moment sequence
$\{a_n\}_{n=0}^\infty$ has a hole. A natural question
arises as to whether the support of a representing
measure $\mu$ of $\{a_n\}_{n=0}^\infty$ has a hole and
what is the relationship between these two holes. We
also study the reverse influence.

If a hole of $\supp\nu$ is properly suited to
$\supp\nu$, then we can locate the corresponding hole
of $\supp\mu$.
   \begin{thm}\label{lem1}
If \eqref{zal} holds and $\nu((\alpha,\beta))=0$ for
some $\alpha, \beta \in \rbb$ such that $0 \Le \alpha
< \beta \Le \gamma:= \sup \supp \nu < \infty$ and
$\alpha \gamma^{\varkappa - 1} < \beta^\varkappa$,
then
   \begin{enumerate}
   \item[(i)] $\mu((\vartheta_1,\vartheta_2))=0$ and
$\vartheta_3=\sup \supp \mu < \infty$,
   \item[(ii)] $\alpha \in \supp \nu$  if and only if
$\vartheta_1 \in \supp \mu$,
   \item[(iii)] $\beta \in \supp \nu$  if and only if
$\vartheta_2 \in \supp \mu$,
   \end{enumerate}
where $\vartheta_1 := \alpha \gamma^{\varkappa - 1}$,
$\vartheta_2:=\beta^\varkappa$ and
$\vartheta_3:=\gamma^\varkappa$.
   \end{thm}
   \begin{proof}
Employing \eqref{compsup} and \eqref{esssup}, we
deduce that the Stieltjes moment sequences
$\{a_n\}_{n=0}^\infty$ and
$\{\sqrtk{a_n}\,\}_{n=0}^\infty$ are determinate,
$\supp \mu \subseteq [0,\gamma^\varkappa]$ and
$\gamma^\varkappa \in \supp \mu$. It is also clear
that $\supp \nu \subseteq [0, \alpha] \sqcup [\beta,
\gamma]$ and $\gamma \in \supp \nu$.

(i) Take $x \in \supp\mu$. By Theorem \ref{suppmu},
there exist $x_1, \ldots, x_\varkappa \in \supp \nu$
such that $x=x_1 \cdots x_\varkappa$. If $x_i \in
[\beta,\gamma]$ for every $i \in \{1, \ldots,
\varkappa\}$, then evidently $x \Ge \beta^\varkappa =
\vartheta_2$. Otherwise, there exists $i_0 \in \{1,
\ldots,\varkappa\}$ such that $x_{i_0} \in
[0,\alpha]$, and thus
   \begin{align*}
x= x_{i_0} \cdot \prod_{j\neq i_0} x_j \Le
\alpha\gamma^{\varkappa-1} = \vartheta_1.
   \end{align*}
Hence $\supp \mu \subseteq [0,\vartheta_1]
\cup[\vartheta_2,\vartheta_3]$, which gives (i).

Now if $\alpha \in \supp \nu$ (respectively:\ $\beta
\in \supp \nu$), then the fact that $\gamma \in \supp
\nu$ enables us to infer from \eqref{porfelTer} that
$\vartheta_1=\alpha \gamma^{\varkappa -1} \in \supp
\mu$ (respectively:\ $\vartheta_2 = \beta^\varkappa\in
\supp \mu$). Therefore, it remains to prove the ``if''
parts of assertions \mbox{(ii)} and~ \mbox{(iii)}.

   (ii) Assume that $\vartheta_1 \in \supp \mu$. Then,
by Theorem \ref{suppmu}, $\vartheta_1=x_1 \cdots
x_\varkappa$ for some $x_1, \ldots, x_\varkappa \in
\supp \nu$. Suppose that, contrary to our claim,
$\alpha \notin \supp \nu$. If $x_i \in [\beta,\gamma]$
for every $i \in \{1, \ldots, \varkappa\}$, then
clearly $\vartheta_1 \Ge \beta^\varkappa=\vartheta_2$,
which contradicts the assumption that $\alpha
\gamma^{\varkappa - 1} < \beta^\varkappa$. Otherwise,
there exists $i_0 \in \{1, \ldots,\varkappa\}$ such
that $x_{i_0} \in [0,\alpha]$. As $\alpha \notin \supp
\nu$, we must have $x_{i_0} < \alpha$ and so
   \begin{align*}
0 < \alpha \gamma^{\varkappa-1} = \vartheta_1 =
x_{i_0} \cdot \prod_{j\neq i_0} x_j < \alpha
\gamma^{\varkappa-1},
   \end{align*}
which is a contradiction.

   (iii) Assume that $\vartheta_2 = \beta^\varkappa
\in \supp \mu$. Suppose that, contrary to our claim,
$\beta \not \in \supp \nu$. Since $\gamma \in \supp
\nu$, we must have $\beta < \gamma$. Clearly, the set
$\supp \nu \cap [\beta,\gamma]$ is compact and
nonempty. Set $\beta^\prime=\min (\supp \nu \cap
[\beta,\gamma])$. Since $\beta \not \in \supp \nu$, we
get $0 \Le \alpha < \beta < \beta^\prime \Le \gamma$
and $\nu((\alpha, \beta^\prime))=0$. Note that
   \begin{align} \label{amk}
\alpha \gamma^{\varkappa-1} < \beta^\varkappa <
\beta^{\prime\varkappa}.
   \end{align}
Hence, by (i), applied to the triplet
$(\alpha,\beta^\prime,\gamma)$, we have $\mu((\alpha
\gamma^{\varkappa-1}, \beta^{\prime\varkappa}))=0$,
which together with \eqref{amk} contradicts the fact
that $\beta^\varkappa \in \supp \mu$.
   \end{proof}
Regarding Theorem \ref{lem1}, one might expect that
the idea of the proof of (iii) would apply to the
proof of (ii) in the case of $\alpha > 0$. However, it
may happen that there is no point in $\supp\nu$ lying
on the left hand side of $\alpha$ (see \eqref{onthel}
in Example \ref{ontheleft}), and so this idea could
not be applied. In turn, (iii) can be proved in the
same manner as (ii). As shown in Example
\ref{UniversalPictures}, the ``if'' parts of
assertions \mbox{(ii)} and \mbox{(iii)} are no longer
true if we drop the assumption that
$\nu((\alpha,\beta))=0$ (though their ``only if''
parts are always true).
   \begin{rem} \label{uwag1}
Under the assumptions of Theorem \ref{lem1}, if
   \begin{align*}
\alpha^\prime :=
   \begin{cases}
\sup([0,\alpha] \cap \sup \nu) & \text{ when }
[0,\alpha] \cap \sup \nu \neq \varnothing,
   \\
0 & \text{ otherwise,}
   \end{cases}
   \end{align*}
and
   \begin{align*}
\beta^\prime := \inf ([\beta,\gamma] \cap \supp \nu),
   \end{align*}
then $0 \Le \alpha^\prime \Le \alpha$, $\beta \Le
\beta^\prime \Le \gamma$,
$\nu((\alpha^\prime,\beta^\prime))=0$ and
$\alpha^\prime \gamma^{\varkappa - 1} <
\beta^{\prime\varkappa}$, which means that the numbers
$\alpha^\prime$, $\beta^\prime$ and $\gamma$ satisfy
the assumptions of Theorem \ref{lem1}.
   \end{rem}
Our goal now is to look for holes of $\supp \nu$ that
may correspond to a given hole of $\supp\mu$.
   \begin{thm}\label{lem1+}
If \eqref{zal} holds and
$\mu((\vartheta_1,\vartheta_2))=0$ for some
$\vartheta_1, \vartheta_2 \in \rbb$ such that $0
\Le \vartheta_1 < \vartheta_2$, then
   \begin{enumerate}
   \item[(i)] $\nu((\beta^\dag, \beta))=0$
provided that $\{a_n\}_{n=0}^\infty$ is
determinate,
   \item[(ii)] $\nu((\alpha,\alpha^\dag)) = 0$
provided that $\vartheta_2 \Le \vartheta_3 :=
\sup \supp \mu < \infty$,
   \item[(iii)] $\nu((\alpha,\beta))=0$
provided that $\vartheta_2 \Le \vartheta_3 := \sup
\supp \mu < \infty$ and any of the following two
conditions holds{\em :}
   \begin{enumerate}
   \item[(iii-a)] either $\beta^\dag < \alpha^\dag$, or
$\beta^\dag \Le \alpha^\dag$ and $\varkappa \Ge 3$, or
$\beta^\dag = \alpha^\dag$ and $\beta \in \supp \nu$,
   \item[(iii-b)]
$\frac{\gamma}{\beta} \, \alpha < \alpha^\dag$ and $\beta
\in \supp \nu$,
   \end{enumerate}
   \end{enumerate}
where $\alpha=\frac{\vartheta_1}{\vartheta_3}
\sqrtk{\vartheta_3}$, $\beta=\sqrtk{\vartheta_2}$,
$\gamma=\sqrtk{\vartheta_3}$, $\alpha^\dag =
\frac{\vartheta_2}{\vartheta_3} \sqrtk{\vartheta_3}$ and
$\beta^\dag =\sqrtk{\vartheta_1}$.
   \end{thm}
   \begin{rem} \label{lem1+uw}
Note that if $\varkappa$ is an integer greater
than or equal to $2$, $0 \Le \vartheta_1 <
\vartheta_2 \Le \vartheta_3 := \sup \supp \mu <
\infty$ and $\beta^\dag \Le \alpha^\dag$, then
$\frac{\gamma}{\beta} \, \alpha < \alpha^\dag$.
Indeed, since the case $\vartheta_2=\vartheta_3$
is obvious, we can assume that $\vartheta_2 <
\vartheta_3$. Then
   \begin{align*}
\frac{\vartheta_1}{\vartheta_3}
 \Le \Big(\frac{\vartheta_2}{\vartheta_3}\Big)^{\varkappa}
< \Big(\frac{\vartheta_2}{\vartheta_3}\Big)^{1+
\frac{1}{\varkappa}},
   \end{align*}
which yields $\frac{\gamma}{\beta} \, \alpha <
\alpha^\dag$. It may happen that condition
\mbox{(iii-b)} of Theorem \ref{lem1+} is
satisfied, while condition \mbox{(iii-a)} does
not hold (cf.\ Example \ref{ontheleft}).
Moreover, assertion \mbox{(iii)} is no longer
true if we drop either the assumption that $\beta
\in \supp \nu$ (cf.\ Example
\ref{UniversalPictures}), or the assumption that
$\frac{\gamma}{\beta} \, \alpha < \alpha^\dag$
(cf.\ Example~ \ref{UniversalPictures2}).
   \end{rem}
   \begin{proof}[Proof of Theorem  \ref{lem1+}]
   (i) Suppose that, contrary to our claim, there
exists $x \in \supp \nu$ such that $ \beta^\dag < x <
\beta$. Then, by \eqref{porfelTer}, we have
$x^\varkappa \in \supp \mu$. This and $\vartheta_1<
x^\varkappa< \vartheta_2$ lead to the contradiction
that $\mu((\vartheta_1,\vartheta_2)) = 0$.

(ii) By \eqref{compsup}, the Stieltjes moment sequence
$\{a_n\}_{n=0}^\infty$ is determinate. Hence, by
Corollary \ref{suppmu-cor}(iii), we have
$\gamma=\sqrtk{\vartheta_3} \in \supp \nu$. Suppose
that, contrary to our claim,
$\nu((\alpha,\alpha^\dag))
> 0$. Then there exists $x\in \supp \nu \cap
(\alpha,\alpha^\dag)$. It follows from
\eqref{porfelTer} that $x \gamma^{\varkappa -1} \in
\supp \mu$. This and the inequalities
   \begin{align*}
\vartheta_1 = \alpha \gamma^{\varkappa -1} < x
\gamma^{\varkappa -1} < \alpha^\dag \gamma^{\varkappa -1} =
\vartheta_2
   \end{align*}
lead to $\mu((\vartheta_1,\vartheta_2)) > 0$, which
contradicts our assumption that
$\mu((\vartheta_1,\vartheta_2)) = 0$.

\mbox{(iii-a)} In view of \eqref{compsup},
$\{a_n\}_{n=0}^\infty$ is determinate. Hence, by (i)
and (ii), we have $\nu((\alpha,\alpha^\dag) \cup
(\beta^\dag,\beta))=0$. So if $\beta^\dag <
\alpha^\dag$, then $\nu((\alpha,\beta))=0$. Assume now
that $\beta^\dag = \alpha^\dag$ and $\varkappa \Ge 3$.
It is enough to show that $\alpha^\dag \notin \supp
\nu$. Suppose that, contrary to our claim,
$\alpha^\dag \in \supp \nu$. By Corollary
\ref{suppmu-cor}(iii), $\gamma = \sqrtk{\vartheta_3}
\in \supp \nu$. This and \eqref{porfelTer} imply that
$(\alpha^\dag)^{\varkappa-1} \gamma \in \supp \mu$.
Since $\varkappa \Ge 3$ and $\alpha^\dag = \beta^\dag
< \beta \Le \gamma$, we have
   \begin{align*}
\vartheta_1 = (\beta^\dag)^{\varkappa} =
(\alpha^\dag)^{\varkappa} <
(\alpha^\dag)^{\varkappa-1} \gamma < \alpha^\dag
\gamma^{\varkappa-1} = \vartheta_2.
   \end{align*}
Hence $\mu((\vartheta_1,\vartheta_2)) > 0$, which is a
contradiction. Finally, the case of $\beta^\dag =
\alpha^\dag$ and $\beta \in \supp \nu$ follows from
Remark \ref{lem1+uw} and \mbox{(iii-b)}.

\mbox{(iii-b)} Set $\alpha_j=(\frac{\gamma}{\beta})^j
\alpha$ and $\alpha_j^\dag=(\frac{\gamma}{\beta})^j
\alpha^\dag$ for $j=0, \dots, \varkappa-1$. Note that
   \begin{align} \label{ajlaj+}
\alpha_j < \alpha_j^\dag \quad \text{for $j=0, \dots,
\varkappa-1$, and} \quad \alpha_0=\alpha < \beta =
\alpha_{\varkappa-1}^\dag.
   \end{align}
We claim that
   \begin{align} \label{bartek}
\nu((\alpha_j,\alpha_j^\dag))=0, \quad j=0, \dots,
\varkappa-1.
   \end{align}
Suppose that, contrary to our claim,
$\nu((\alpha_j,\alpha_j^\dag)) > 0$ for some $j\in
\{0, \dots, \varkappa-1\}$. Then there exists $x \in
(\alpha_j,\alpha_j^\dag) \cap \supp \nu$. By Corollary
\ref{suppmu-cor}(iii), $\gamma \in \supp \nu$. Since
$\beta, \gamma \in \supp \nu$, we infer from
\eqref{porfelTer} that $x \beta^j \gamma^{\varkappa -
1 -j} \in \supp \mu$. Noting that
   \begin{multline*}
\vartheta_1 = \alpha \gamma^{\varkappa-1} =
\Big(\frac{\gamma}{\beta}\Big)^j \alpha \beta^j
\gamma^{\varkappa -1 - j} = \alpha_j \beta^j
\gamma^{\varkappa -1 - j}
   \\
   < x \beta^j \gamma^{\varkappa - 1 -j} < \alpha_j^\dag
\beta^j \gamma^{\varkappa -1 - j} =
\Big(\frac{\gamma}{\beta}\Big)^j \alpha^\dag \beta^j
\gamma^{\varkappa -1 - j} = \alpha^\dag
\gamma^{\varkappa-1} = \vartheta_2,
   \end{multline*}
we get $\mu((\vartheta_1, \vartheta_2)) > 0$, in
contradiction with our assumption that $\mu((\vartheta_1,
\vartheta_2)) = 0$.

The inequality $\frac{\gamma}{\beta} \, \alpha <
\alpha^\dag$ is easily seen to be equivalent to
   \begin{align} \label{eqdag}
\alpha_{j+1} < \alpha_j^\dag, \quad j=0, \ldots,
\varkappa-2.
   \end{align}
We will show that
   \begin{align} \label{sum}
(\alpha, \beta) = \bigcup_{j=0}^{\varkappa-1}
(\alpha_j,\alpha_j^\dag).
   \end{align}
Indeed, assuming that\footnote{\;The case of
$\vartheta_2 = \vartheta_3$ is obvious.} $\vartheta_2 <
\vartheta_3$, we see that if $x \in
(\alpha,\alpha_{\varkappa-1}]$, then $x \in (\alpha_j,
\alpha_{j+1}]$ for some $j \in \{0, \ldots,
\varkappa-2\}$. Hence, by \eqref{eqdag}, $x \in
(\alpha_j, \alpha_{j}^\dag)$. Since, by \eqref{ajlaj+},
$(\alpha_{\varkappa-1},\alpha_{\varkappa-1}^\dag) =
(\alpha_{\varkappa-1},\beta)$, the equality \eqref{sum}
is proved.

Combining \eqref{bartek} with \eqref{sum}, we conclude that
$\nu((\alpha,\beta))=0$. This completes the proof.
   \end{proof}
Regarding Theorem \ref{lem1+}, we can write the
following analogue of Remark \ref{uwag1}.
   \begin{rem}
Suppose that \eqref{zal} holds and $\vartheta_1,
\vartheta_2, \vartheta_3$ are real numbers such that
$0 \Le \vartheta_1 < \vartheta_2 \Le \vartheta_3 =
\sup \supp \mu$ and
$\mu((\vartheta_1,\vartheta_2))=0$. Set
$\vartheta_1^\prime = \vartheta_1$,
$\vartheta_2^\prime = \inf([\vartheta_2,\vartheta_3]
\cap \supp \mu)$ and $\vartheta_3^\prime =
\vartheta_3$. Then $0 \Le \vartheta_1^\prime <
\vartheta_2^\prime \Le \vartheta_3^\prime = \sup \supp
\mu$ and $\mu((\vartheta_1^\prime,
\vartheta_2^\prime))=0$. Let $\alpha, \alpha^\dag,
\beta, \beta^\dag, \gamma$ (respectively
$\alpha^\prime, \alpha^{\dag\prime}, \beta^\prime,
\beta^{\dag\prime}, \gamma^\prime$) be numbers
attached to $\vartheta_1, \vartheta_2, \vartheta_3$
(respectively $\vartheta_1^\prime, \vartheta_2^\prime,
\vartheta_3^\prime$) via \eqref{theta2.5} and
\eqref{abplus}. Then $\alpha = \alpha^\prime$,
$\alpha^\dag \Le \alpha^{\dag\prime}$, $\beta \Le
\beta^\prime$, $\beta^{\dag} = \beta^{\dag\prime}$ and
$\gamma = \gamma^\prime$. Hence, if $\beta^\dag <
\alpha^\dag$ or $\beta^\dag \Le \alpha^\dag$, then
$\beta^\dag < \alpha^{\dag\prime}$ or $\beta^\dag \Le
\alpha^{\dag\prime}$ respectively. This means that we
can apply Theorem \ref{lem1+} to the new system of
numbers $\vartheta_1, \vartheta_2^\prime, \vartheta_3,
\alpha, \alpha^{\dag\prime}, \beta^\prime, \beta^\dag,
\gamma$ (note that if $\beta \in \supp \nu$, then by
Corollary \ref{suppmu-cor}\,(i) we have $\vartheta_2 =
\beta^\varkappa \in \supp \mu$, and so
$\vartheta_2^\prime=\vartheta_2$).
   \end{rem}
   \begin{cor} \label{cor1}
Suppose \eqref{zal} holds. Let $\vartheta_1,
\vartheta_2, \vartheta_3$ be real numbers such that $0
\Le \vartheta_1 < \vartheta_2 \Le \vartheta_3$ and let
$\alpha:=\frac{\vartheta_1}{\vartheta_3}
\sqrtk{\vartheta_3}$, $\beta:=\sqrtk{\vartheta_2}$,
$\gamma:=\sqrtk{\vartheta_3}$ and $\beta^\dag
:=\sqrtk{\vartheta_1}$. Then the following assertions
are valid.
   \begin{enumerate}
   \item[(i)] If $\vartheta_2
< \vartheta_3$, $\{a_n\}_{n=0}^\infty$ is determinate
and $\mu((\vartheta_1,\vartheta_2))=
\mu((\vartheta_2,\vartheta_3)) = 0$, then
$\nu((\beta^\dag, \beta)) = \nu((\beta,\gamma)) = 0$.
   \item[(ii)] If $\vartheta_2=\vartheta_3$,
$\vartheta_2 = \sup \supp \mu$ and
$\mu((\vartheta_1,\vartheta_2)) = 0$, then
$\nu((\alpha,\beta)) = 0$ and $\beta \in \supp \nu$.
   \item[(iii)] If $\vartheta_2=\vartheta_3$,
$\vartheta_1 \in \supp \mu$, $\vartheta_2 = \sup \supp
\mu$ and $\mu((\vartheta_1,\vartheta_2)) = 0$, then
$\nu((\alpha,\beta)) = 0$ and $\alpha, \beta \in \supp
\nu$.
   \item[(iv)] If $\vartheta_2=\vartheta_3$,
then the following two conditions are equivalent\/{\em
:}
   \begin{enumerate}
   \item[(a)] $\vartheta_1 \in \supp \mu$,
$\vartheta_2 = \sup \supp \mu$ and
$\mu((\vartheta_1,\vartheta_2))=0$,
   \item[(b)]
$\alpha \in \supp \nu$, $\beta=\sup \supp \nu$ and
$\nu((\alpha,\beta))=0$.
   \end{enumerate}
   \end{enumerate}
   \end{cor}
   \begin{proof}
(i) Apply Theorem \ref{lem1+}\,(i) to intervals
$(\vartheta_1,\vartheta_2)$ and
$(\vartheta_2,\vartheta_3)$.

(ii) and (iii) These two assertions can be deduced from
assertion (ii) of Theorem \ref{lem1+}, \eqref{esssup}
and assertions (ii) and (iii) of Theorem \ref{lem1}
(because $\beta = \alpha^\dag$).

   (iv) Implication (a)$\Rightarrow$(b) can be inferred
from (iii) and \eqref{esssup}. The reverse implication
is an immediate consequence of Theorem \ref{lem1}.
   \end{proof}
See Appendix for another proof of the implication
(a)$\Rightarrow$(b) of Corollary \ref{cor1}\,(iv).

We conclude this section by discussing the case when
the support of a representing measure of a Stieltjes
moment sequence is contained in a closed interval
$[\vartheta,\infty)$, $\vartheta > 0$. Note that
$\vartheta$ corresponds to $\vartheta_2$ in Theorem
   \ref{lem1}.
   \begin{pro}
Suppose that {\em \eqref{zal}} holds and
$\{a_n\}_{n=0}^\infty$ is determinate. Then the
following assertions are valid.
   \begin{enumerate}
   \item[(i)] If $\nu([0,\beta))=0$ for some
real number $\beta > 0$, then
$\mu([0,\beta^\varkappa))=0$.
   \item[(ii)] If $\mu([0,\vartheta))=0$ for some
real number $\vartheta > 0$, then
$\nu([0,\sqrtk{\vartheta}))=0$.
   \item[(iii)] If $\mu([0,\vartheta))=0$ for some
real number $\vartheta > 0$ such that $\vartheta \in
\supp \mu$, then $\nu([0,\sqrtk{\vartheta}))=0$ and
$\sqrtk{\vartheta} \in \supp \nu$.
   \end{enumerate}
   \end{pro}
   \begin{proof}
   Assertions (i) and (ii) follow from Theorem
\ref{suppmu}.

   (iii) By (ii), we have
$\nu([0,\sqrtk{\vartheta}))=0$, and thus $0 \notin
\supp\nu$. It follows from Theorem \ref{suppmu} that
$\vartheta = x_1 \cdots x_\varkappa$ for some $x_1,
\ldots, x_\varkappa \in \supp \nu$. Suppose that,
contrary to our claim, $\sqrtk{\vartheta} \notin \supp
\nu$. Then $x_j > \sqrtk{\vartheta}$ for all $j\in \{1,
\ldots, \varkappa\}$, and thus $\vartheta = x_1 \cdots
x_\varkappa > \vartheta$, which is a contradiction.
   \end{proof}
   \section{When is $(\alpha, \beta)$ a hole of
$\supp\nu$?\label{secnew}}
   Before stating a theorem which provides an answer
to the above question, we introduce a new parameter
$\iota_{\mathrm{s}}^*$ and prove two technical lemmas
about the parameters $\iota_{\mathrm{s}}$ (cf.\
\eqref{iotadef}) and $\iota_{\mathrm{s}}^*$. For real
numbers $\vartheta_1, \vartheta_2, \vartheta_3$ such
that $0 < \vartheta_1 < \vartheta_2 < \vartheta_3$,
we~ set
   \begin{align} \label{iota*}
\iota_{\mathrm{s}}^* =
\iota_{\mathrm{s}}^*(\vartheta_1, \vartheta_2,
\vartheta_3) := 1 + \left\lfloor
\frac{\log{(\vartheta_3/\vartheta_2)}}
{\log{(\vartheta_2/\vartheta_1)}}\right\rfloor \Ge 1.
   \end{align}
Note that the following equalities hold for all $t \in
(0,\infty)$,
   \begin{align*}
\text{$\iota_{\mathrm{s}}(\vartheta_1, \vartheta_2,
\vartheta_3) = \iota_{\mathrm{s}}\big(t \vartheta_1, t
\vartheta_2, t \vartheta_3\big)$ and
$\iota_{\mathrm{s}}^*(\vartheta_1, \vartheta_2,
\vartheta_3) = \iota_{\mathrm{s}}^* \big(t
\vartheta_1, t \vartheta_2, t \vartheta_3\big)$.}
   \end{align*}
   \begin{lem} \label{techn}
If $\vartheta_1, \vartheta_2, \vartheta_3 \in \rbb$ are
such that $0 < \vartheta_1 < \vartheta_2 <
\vartheta_3$, then the following assertions are
valid{\em :}
   \begin{enumerate}
   \item[(i)] $\iota_{\mathrm{s}} = 2$ implies
$\iota_{\mathrm{s}}^* \Ge 2$,
   \item[(ii)] $\iota_{\mathrm{s}} = 3$ implies
$\iota_{\mathrm{s}}^* \in \{1,2\}$,
   \item[(iii)] $\iota_{\mathrm{s}} = 3$ and
$\iota_{\mathrm{s}}^*=2$ if and only if
$\vartheta_1\vartheta_3 = \vartheta_2^2$,
   \item[(iv)] $\iota_{\mathrm{s}} \Ge 4$ implies
$\iota_{\mathrm{s}}^*=1$,
   \item[(v)] $\iota_{\mathrm{s}}^* =1$ if and
only if $\vartheta_1\vartheta_3 < \vartheta_2^2$,
   \item[(vi)] $\iota_{\mathrm{s}}^* =1$ implies
$\iota_{\mathrm{s}} \Ge 3$,
   \item[(vii)] $\iota_{\mathrm{s}}^* =2$ implies
$\iota_{\mathrm{s}} \in \{2,3\}$,
   \item[(viii)] $\iota_{\mathrm{s}}^* \Ge 3$ implies
$\iota_{\mathrm{s}} = 2$.
   \end{enumerate}
Moreover, for every integer $p \Ge 2$, there exist
real numbers $\vartheta_1, \vartheta_2, \vartheta_3$
such that $0 < \vartheta_1 < \vartheta_2 <
\vartheta_3$, $\iota_{\mathrm{s}}=2$ and
$\iota_{\mathrm{s}}^*=p$.
   \end{lem}
   \begin{proof}
Suppose that $\vartheta_1, \vartheta_2, \vartheta_3
\in \rbb$ are such that $0 < \vartheta_1 < \vartheta_2
< \vartheta_3$. First we show that there exists $r \in
[0,1)$ such that
   \begin{align}  \label{iota+iota*}
\iota_{\mathrm{s}} -2 + r > 0 \text{ and }
\iota_{\mathrm{s}}^* = 1 + \left\lfloor
\frac{1}{\iota_{\mathrm{s}} -2 + r}\right\rfloor.
   \end{align}
Indeed, by definition of $\iota_{\mathrm{s}}$, there
exists $r \in [0,1)$ such that $\iota_{\mathrm{s}} -2
+ r > 0$ and $\frac{\log{(\vartheta_3/\vartheta_1)}}
{\log{(\vartheta_3/\vartheta_2)}} = \iota_{\mathrm{s}}
-1+r$. This gives us $\frac{\vartheta_3} {\vartheta_1}
= \big(\frac{\vartheta_3}
{\vartheta_2}\big)^{\iota_{\mathrm{s}}-1+r}$. As a
consequence, we have
   \begin{align} \label{iota+iota*2}
\frac{\vartheta_3}{\vartheta_2} =
\Big(\frac{\vartheta_2}{\vartheta_1}\Big)^{\frac{1}
{\iota_{\mathrm{s}} -2 + r}},
   \end{align}
which yields \eqref{iota+iota*}.

Conditions (i) to (viii) except for (iii) and (v) can
be deduced from \eqref{iota+iota*}. Condition (iii)
follows from \eqref{iota+iota*} and
\eqref{iota+iota*2}. Clearly, $\iota_{\mathrm{s}}^* =1$
if and only if $\frac{\log{(\vartheta_3/\vartheta_2)}}
{\log{(\vartheta_2/\vartheta_1)}} < 1$, or equivalently
if and only if $\vartheta_1\vartheta_3 <
\vartheta_2^2$, which gives (v).

Now we justify the ``moreover'' part. Set
$\vartheta_3=1$. For fixed $\vartheta_1 \in (0, 1)$
and $p\in \{2,3, \ldots\}$, we define $\vartheta_2 =
\vartheta_1^{\frac{2p-1}{2p+1}}$. It is a simple
matter to verify that $\vartheta_1 < \vartheta_2 <
\vartheta_3$, $\iota_{\mathrm{s}}=2$ and
$\iota_{\mathrm{s}}^*=p$. This completes the proof.
   \end{proof}
   \begin{lem} \label{techn+}
Let $\vartheta_1, \vartheta_2, \vartheta_3$ be
real numbers such that $0 < \vartheta_1 <
\vartheta_2 < \vartheta_3$. Then the following
assertions hold\/{\em :}
   \begin{enumerate}
   \item[(i)] if $\alpha^\dag = \beta^\dag$ and
$\varkappa \Ge 3$, then $\iota_{\mathrm{s}}^*
=1$,
   \item[(ii)] if  $\varkappa=2$, then $\alpha^\dag =
\beta^\dag$ if and only if $\iota_{\mathrm{s}} =
3$ and $\iota_{\mathrm{s}}^* = 2$,
   \item[(iii)] if $\varkappa \Ge
\iota_{\mathrm{s}}^*$, then $\frac{\gamma}{\beta} \,
\alpha < \alpha^\dag$,
   \end{enumerate}
where $\alpha$, $\alpha^\dag$, $\beta$, $\beta^\dag$
and $\gamma$ are as in \eqref{theta2.5} and
\eqref{abplus}.
   \end{lem}
   \begin{proof}
(i) By \eqref{abplus}, we have
   \begin{align} \label{eqthet}
\vartheta_2 = \vartheta_3
\sqrtk{\frac{\vartheta_1}{\vartheta_3}}.
   \end{align}
Since $\varkappa \Ge 3$, we get $1-\frac{2}{\varkappa}
> 0$, and thus $\vartheta_1^{1-\frac{2}{\varkappa}} <
\vartheta_3^{1-\frac{2}{\varkappa}}$, which together
with \eqref{eqthet} implies that $\vartheta_1
\vartheta_3 < \Big(\vartheta_3
\sqrtk{\frac{\vartheta_1}{\vartheta_3}}\,\Big)^2 =
\vartheta_2^2$. By Lemma \ref{techn}(v),
$\iota_{\mathrm{s}}^* =1$.

(ii) By \eqref{abplus}, $\alpha^\dag =
\beta^\dag$ if and only if
$\vartheta_1\vartheta_3 = \vartheta_2^2$.
Applying Lemma \ref{techn}(iii) completes the
proof of (ii).

(iii) Since $\varkappa \Ge \iota_{\mathrm{s}}^*$,
we see that $\varkappa >
\frac{\log{(\vartheta_3/\vartheta_2)}}
{\log{(\vartheta_2/\vartheta_1)}}$. Hence
$\frac{\vartheta_2} {\vartheta_1}
>\big(\frac{\vartheta_3} {\vartheta_2}\big)^{\frac{1}
{\varkappa}}$, and thus
$\frac{\vartheta_3}{\vartheta_1}
> \big(\frac{\vartheta_3}{\vartheta_2}\big)^{1 +
\frac{1}{\varkappa}}$, which implies that
$\frac{\gamma}{\beta} \, \alpha < \alpha^\dag$.
   \end{proof}
Theorem \ref{17.02.11} below gives sufficient
conditions for an interval $(\alpha, \beta)$ to be a
hole of $\supp \nu$ written in terms of the parameters
$\iota_{\mathrm{s}}$ and $\iota_{\mathrm{s}}^*$.
   \begin{thm} \label{17.02.11}
If \eqref{zal} holds, $\mu((\vartheta_1,
\vartheta_2))=0$ for some $\vartheta_1, \vartheta_2
\in \rbb$ such that $0 < \vartheta_1 < \vartheta_2 <
\vartheta_3:=\sup \supp \mu < \infty$ and any of the
conditions {\em (i)}-{\em (v)} below is satisfied,
then $\nu((\alpha,\beta))=0$.
   \begin{enumerate}
   \item[(i)] $\varkappa \Ge \iota_{\mathrm{s}}^*$
and $\beta \in \supp \nu$.
   \item[(ii)] $\iota_{\mathrm{s}} \Ge
\iota_{\mathrm{s}}^*$ and $\beta \in \supp \nu$.
   \item[(iii)] $\iota_{\mathrm{s}} \Ge 3$ and
$\beta \in \supp \nu$.
   \item[(iv)] $\iota_{\mathrm{s}} \Ge 4$.
   \item[(v)] $\iota_{\mathrm{s}}^*=1$.
   \end{enumerate}
$($$\alpha, \beta, \iota_{\mathrm{s}},
\iota_{\mathrm{s}}^*$ are as in \eqref{theta2.5},
\eqref{iotadef} and \eqref{iota*}.$)$
   \end{thm}
   \begin{proof}
(i) Apply Lemma \ref{techn+}(iii) and Theorem
\ref{lem1+}\,\mbox{(iii-b)}.

(v) We deal first with the case in which $\vartheta_2
\in \supp\mu$. We will show that $\beta \in \supp\nu$.
By \eqref{compsup}, Corollary \ref{suppmu-cor}\,(iii)
and Theorem \ref{suppmu}, there exist $x_1, \ldots,
x_\varkappa \in \supp\nu$ such that $\vartheta_2=x_1
\cdots x_\varkappa$. It follows from \eqref{abplus}
and \eqref{duzoal} that $\vartheta_2 = \alpha^\dag
\gamma^{\varkappa-1}$ and $0<\alpha^\dag, x_1, \ldots,
x_\varkappa \Le \gamma$. Hence
   \begin{align} \label{abc1}
x_j \in [\alpha^\dag,\gamma], \quad j \in \{1, \ldots,
\varkappa\}.
   \end{align}
Note that $x_1= \cdots = x_\varkappa$. Indeed,
otherwise $x_k < x_l$ for some $k\neq l$. Since, by
Lemma \ref{techn}(v), $\vartheta_1 \vartheta_3 <
\vartheta_2^2$, we have
   \begin{align}  \label{Jurek}
\vartheta_2 > x_k^2 \cdot \prod_{j \notin\{k,l\}} x_j
= \vartheta_2 \frac{x_k}{x_l}
\overset{\eqref{abc1}}\Ge \vartheta_2
\frac{\alpha^\dag}{\gamma} =
\frac{\vartheta_2^2}{\vartheta_3} > \vartheta_1.
   \end{align}
Applying Theorem \ref{suppmu} we see that $x_k^2 \cdot
\prod_{j \notin\{k,l\}} x_j \in \supp\mu$, which
together with \eqref{Jurek} shows that
$\mu((\vartheta_1, \vartheta_2))
> 0$. This leads to a contradiction. Since $x_1=
\cdots = x_\varkappa$, we deduce that $\beta =
\sqrtx{1}{\varkappa}{\vartheta_2} = x_1 \in \supp\nu$.
Applying (i), we get $\nu((\alpha,\beta))=0$.

If $\vartheta_2 \notin \supp\mu$, then we argue as
follows. Set $\vartheta_2^\prime =
\inf([\vartheta_2,\vartheta_3] \cap \supp \mu)$
and $\beta^\prime=\sqrtk{\vartheta_2^\prime}$.
Then $\beta \Le \beta^\prime$ and
$\mu((\vartheta_1,\vartheta_2^\prime))=0$. If
$\vartheta_2^\prime = \vartheta_3$, then by
Corollary \ref{cor1}(ii) we have
$\nu((\alpha,\beta)) \Le
\nu((\alpha,\beta^\prime))=0$. Otherwise
$\vartheta_2^\prime < \vartheta_3$. Since $1 \Le
\iota_{\mathrm{s}}^*(\vartheta_1,
\vartheta_2^\prime, \vartheta_3) \Le
\iota_{\mathrm{s}}^*(\vartheta_1, \vartheta_2,
\vartheta_3)=1$ and $\vartheta_2^\prime \in \supp
\mu$, we may apply the argument from the previous
paragraph to $\vartheta_1, \vartheta_2^\prime,
\vartheta_3$. What we obtain is
$\nu((\alpha,\beta)) \Le
\nu((\alpha,\beta^\prime))=0$.

(iv) Apply (v) and Lemma \ref{techn}(iv).

(ii) If $\varkappa \Ge \iota_{\mathrm{s}}^*$, then we
may apply (i). Consider the case of $\varkappa <
\iota_{\mathrm{s}}^*$. Then, by $\iota_{\mathrm{s}} \Ge
\iota_{\mathrm{s}}^*$, we have $2 \Le \varkappa <
\iota_{\mathrm{s}}$. Thus, if $\iota_{\mathrm{s}} \Ge
4$, then we may apply (iv). It remains to consider the
case of $\iota_{\mathrm{s}} = 3$. Then clearly
$\varkappa = 2$. It follows from $\iota_{\mathrm{s}} =
3$ that $\frac{\log{(\vartheta_3/\vartheta_1)}}
{\log{(\vartheta_3/\vartheta_2)}} \Ge 2$, or
equivalently that $\frac{\vartheta_3}{\vartheta_1} \Ge
(\frac{\vartheta_3}{\vartheta_2})^2$, which is
equivalent to $\beta^\dag \Le \alpha^\dag$. Applying
Theorem \ref{lem1+}\,\mbox{(iii-a)} gives
$\nu((\alpha,\beta))=0$.

(iii) Since, by parts (ii) and (iv) of Lemma
\ref{techn}, the inequality $\iota_{\mathrm{s}} \Ge 3$
implies that $\iota_{\mathrm{s}} >
\iota_{\mathrm{s}}^*$, we can apply (ii). This
completes the proof.
   \end{proof}
Note that conditions (ii)-(v) of Theorem
\ref{17.02.11} impose no restriction on $\varkappa$,
and that Theorem \ref{17.02.11} is no longer true if
any of the assumptions (i)-(v) is dropped (see Example
\ref{UniversalPictures} for the discussion concerning
the assumption $\beta \in \supp \nu$, and Example
\ref{UniversalPictures2} for the discussion concerning
the remaining assumptions).
   \begin{cor} \label{co2}
Let $\{a_n\}_{n=0}^\infty$ be a Stieltjes moment
sequence with a representing measure $\mu$ such that
$0 < \vartheta_3:=\sup \supp \mu < \infty$. Suppose
$\mu((\vartheta_1, \vartheta_2))=0$ for some
$\vartheta_1, \vartheta_2 \in \rbb$ such that $0 <
\vartheta_1 < \vartheta_2 < \vartheta_3$. Assume that
the set
   \begin{align} \label{setJ}
J := \Big\{\varkappa \in \zbb_+\colon \varkappa \Ge 2
\text{ and } \big\{\sqrtx{1}{\varkappa}
{a_n}\,\big\}_{n=0}^\infty \text{ is a Stieltjes
moment sequence}\Big\}
   \end{align}
is nonempty. If $\iota_{\mathrm{s}}^* = 1$, then the
following conditions are equivalent\/{\em :}
   \begin{enumerate}
   \item[(i)] $\vartheta_2 \in \supp \mu$,
   \item[(ii)] $\beta(\varkappa) \in \supp \nu_\varkappa$
for some $\varkappa \in J$,
   \item[(iii)] $\beta(\varkappa) \in \supp \nu_\varkappa$
for every $\varkappa \in J$,
   \end{enumerate}
where $\beta(\varkappa)=\sqrtk{\vartheta_2}$ and
$\nu_\varkappa$ is a representing measure of
$\{\sqrtx{1}{\varkappa} {a_n}\}_{n=0}^\infty$.
   \end{cor}
   \begin{proof}
By \eqref{esssup}, we have $\gamma(\varkappa):=
\sqrtk{\vartheta_3} = \sup \supp \nu_{\varkappa}$ for
every $\varkappa \in J$.

(i)$\Rightarrow$(iii) It follows form Theorem
\ref{17.02.11}\,(v) that
$\nu_{\varkappa}((\alpha(\varkappa),\beta(\varkappa)))=0$
for every $\varkappa \in J$, where $\alpha(\varkappa) =
\frac{\vartheta_1}{\vartheta_3} \sqrtk{\vartheta_3}$.
Hence, by Theorem \ref{lem1}\,(iii), $\beta(\varkappa)
\in \supp \nu_\varkappa$ for every $\varkappa \in J$.

(iii)$\Rightarrow$(ii) Evident.

(ii)$\Rightarrow$(i) Applying \eqref{porfelTer}, we
see that $\vartheta_2 = \beta(\varkappa)^\varkappa \in
\supp \mu$.
   \end{proof}
It is worth mentioning that implication
(i)$\Rightarrow$(iii) of Corollary \ref{co2} is no
longer true if we drop the assumption that
$\iota_{\mathrm{s}}^* = 1$ (cf.\ Example
\ref{UniversalPictures}), though the reverse
implication is always true. What is more, the set $J$
defined in \eqref{setJ} may not be a set of
consecutive integers (cf.\ Example \ref{prznotsq}). As
shown in Example \ref{UniversalPictures}, the set $J$
may contain only one point. It may also happen that
$J=\{2,3,4, \ldots\}$ (cf.\ \cite{Hor2}).
   \section{\label{przyk}Examples}
In this section we gather examples that illustrate the
delicate nature of results appearing in Sections
\ref{sec3}, \ref{sec4} and \ref{secnew}. In what follows we
adhere to the notation in \eqref{theta2.5}, \eqref{abplus},
\eqref{iotadef} and \eqref{iota*}. If
$\{a_n\}_{n=0}^\infty$ and $\{\sqrtk{a_n}\,\}_{n=0}^\infty$
are determinate Stieltjes moment sequences, then their
representing measures will be denoted by $\mu$ and $\nu$
respectively.

We begin by showing that the closure sign in
\eqref{porfelTer} cannot be omitted.
   \begin{exa} \label{notclos}
Let $\varkappa=2$. For $\tau \in \{0, \frac12\}$, we
set
   \begin{align} \label{exp}
\sqrt{a_n} = \sum_{j=2}^\infty \frac{1}{2^j(j+\tau)^n}
+ \sum_{j=2}^\infty \frac{j^n}{\E^{j^2}}, \quad n \in
\zbb_+.
   \end{align}
It is clear that $\{\sqrt{a_n}\,\}_{n=0}^\infty$ and
$\{a_n\}_{n=0}^\infty$ are Stieltjes moment sequences,
and that $\nu := \sum_{j=2}^\infty 2^{-j}
\delta_{\frac{1}{j+\tau}} + \sum_{j=2}^\infty \E^{-j^2}
\delta_{j}$ is a representing measure of
$\{\sqrt{a_n}\,\}_{n=0}^\infty$. Hence
   \begin{align} \label{supni}
\supp \nu = \{0\} \cup \Big\{\ldots, \frac{1}{4+\tau},
\frac{1}{3+\tau}, \frac{1}{2+\tau}\Big\} \cup \{2,3,4,
\ldots\}.
   \end{align}
Now we show that $\{a_n\}_{n=0}^\infty$ is a
determinate Stieltjes moment sequence. For $n\Ge 1$,
we define the function $f_{n}\colon \rbb_+ \to \rbb_+$
by $f_n(x)=x^{2n}\E^{-x^2}$ for $x\in \rbb_+$. It is
easily seen that $f_n$ is strictly increasing on the
interval $[0,\sqrt{n}\,]$ and strictly decreasing on
the interval $[\sqrt{n},\infty)$. This implies that
   \begin{align*}
\sqrt{a_{2n}} & \overset{\eqref{exp}} \Le 1 +
\sum_{j=2}^\infty f_n(j)
   \\
& \hspace{1ex}\Le 1 +
\sum_{j=1}^{\lfloor\sqrt{n}\rfloor-1} f_n(j) +
\sum_{j=\lfloor\sqrt{n}\rfloor+ 2}^\infty f_n(j) +
f_n(\lfloor\sqrt{n}\rfloor) +
f_n(\lfloor\sqrt{n}\rfloor+1)
   \\
& \hspace{1ex} \Le 1 + \int_0^\infty f_n(x) \,\D x + 2
f_n(\sqrt{n}), \quad n \Ge 4.
   \end{align*}
By a suitable change of variables, we have
   \begin{align*}
\int_0^\infty f_n(x) \,\D x & = \frac 12 \int_0^\infty
x^{n-\frac 12} \E^{-x} \,\D x \Le \int_0^1 x^{n-\frac
12} \E^{-x} \,\D x + \int_1^\infty x^{n-\frac 12}
\E^{-x} \,\D x
   \\
& \Le 1 + \int_1^\infty x^{n} \E^{-x} \,\D x \Le 1 +
\int_0^\infty x^{n} \E^{-x} \,\D x = 1+n!, \quad n \Ge
1.
   \end{align*}
It is also clear that $f_n(\sqrt{n}) \Le n^n$ for all
$n \Ge 1$. Putting all these together yields
$\sqrt{a_{2n}} \Le 4n^n$ for all $n \Ge 4$, which
implies that $\sum_{n=1}^\infty
a_{2n}^{-\nicefrac{1}{2n}} = \infty$. Hence, by the
Carleman criterion (cf.\ \cite[Theorem 1.10]{sh-tam}),
the Stieltjes moment sequence $\{a_n\}_{n=0}^\infty$ is
determinate.

We first consider the case of $\tau=0$. By \eqref{supni},
the set $\pi_\varkappa\big((\supp \nu)^\varkappa\big)$
coincides with the set of all nonnegative rational numbers,
and so, in view of Theorem \ref{suppmu}, we have $\supp \mu
= \rbb_+$.

Now suppose that $\tau=\nicefrac 12$. It follows from
\eqref{supni} that $1 \notin \pi_\varkappa\big((\supp
\nu)^\varkappa\big)$. Since $1=\lim_{j\to\infty}
\frac{j}{j+\frac 12}$, we infer from \eqref{porfelTer}
that $1 \in \supp \mu$.

As above we verify that the Stieltjes moment sequence
$\{a_n\}_{n=0}^\infty$ given by
   \begin{align*}
\sqrt{a_n} = \sum_{j=2}^\infty \frac{1}{2^j(j+\tau)^n}
+ \sum_{j=2}^\infty \frac{j^n}{\E^{j}}, \quad n \in
\zbb_+ \quad \Big(\tau \in \Big\{0, \frac 12
\Big\}\Big),
   \end{align*}
satisfies the inequality $\sqrt{a_{n}} \Le 3n^n$ for
all $n \Ge 2$, which implies that $\sum_{n=1}^\infty
a_{n}^{-\nicefrac{1}{2n}} = \infty$. Hence, by the
Carleman criterion (cf.\ \cite[Theorem 1.11]{sh-tam}),
the Stieltjes moment sequence $\{a_n\}_{n=0}^\infty$
is determinate (in fact, it is determinate as a
Hamburger moment sequence, cf.\ \cite[Corollary
4.5]{sim}). Choosing $\tau \in \{0, \frac 12\}$, we
obtain new examples of Stieltjes moment sequences with
the required properties.
   \end{exa}
The next example is related to Corollary
\ref{suppmu-cor} and Theorems \ref{lem1} and
\ref{lem1+}.
   \begin{exa} \label{ontheleft}
Set $\varkappa=3$. Let $\vartheta_1, \vartheta_2,
\vartheta_3$ be real numbers such that
   \begin{align*}
0 < \vartheta_1 < \vartheta_2 < \vartheta_3, \quad
\frac{\vartheta_2}{\vartheta_3} <
\sqrt[3]{\frac{\vartheta_1}{\vartheta_3}} \quad
\text{and} \quad
\sqrt{\frac{\vartheta_1}{\vartheta_3}} <
\sqrt[3]{\frac{\vartheta_2}{\vartheta_3}}
   \end{align*}
(e.g., $\vartheta_1 = \frac{1}{\sqrt[3]{3}}$,
$\vartheta_2=1$ and $\vartheta_3 = 2$). Then $0 <
\alpha < \alpha^\dag < \beta^\dag < \beta < \gamma$
and
   \begin{align}      \label{nier}
   \begin{aligned} \alpha^3 < \alpha^2 \beta <
\alpha^2\gamma < \alpha \beta^2 & < \alpha \beta \gamma
   \\
& < \alpha \gamma^2=\vartheta_1 < \vartheta_2=\beta^3 <
\beta^2 \gamma < \beta\gamma^2 < \gamma^3=\vartheta_3.
   \end{aligned}
   \end{align}
Set $a_n = (\alpha^n + \beta^n + \gamma^n)^3$ for $n
\in \zbb_+$. Clearly $\{a_n\}_{n=0}^\infty$ and
$\{\sqrt[3]{a_n}\}_{n=0}^\infty$ are determinate
Stieltjes moment sequences. We easily verify that the
terms of the sequence \eqref{nier} form the support of
$\mu$ and that $\supp \nu = \{\alpha, \beta,
\gamma\}$. Thus $\mu((\vartheta_1,\vartheta_2))=0$,
$\{\vartheta_1,\vartheta_2\} \subseteq \supp \mu$ and
$\mu([0,\vartheta_1))>0$, though
   \begin{align} \label{onthel}
\text{$\nu((\alpha,\beta))=0$, $\{\alpha,\beta\}
\subseteq \supp \nu$ and $\nu([0,\alpha))=0$.}
   \end{align}
Since $\card(\supp \mu) = 10$, we see that $\supp \nu
\varsubsetneq \sqrtk{\supp \mu}$. Moreover, if we
replace the inequality
$\sqrt{\frac{\vartheta_1}{\vartheta_3}} <
\sqrt[3]{\frac{\vartheta_2}{\vartheta_3}}$ by the
stronger one
$\sqrt[4]{\frac{\vartheta_1}{\vartheta_3}} <
\sqrt[3]{\frac{\vartheta_2}{\vartheta_3}}$ (which is
still satisfied by $\vartheta_1 =
\frac{1}{\sqrt[3]{3}}$, $\vartheta_2=1$ and
$\vartheta_3 = 2$), then $\frac{\gamma}{\beta} \alpha
< \alpha^\dag$ and $\beta \in \supp \nu$, though
$\beta^\dag > \alpha^\dag$.
   \end{exa}
Example \ref{UniversalPictures} below shows that the
assumption that one or two endpoints of the interval
$(\alpha,\beta)$ belong to $\supp\nu$ is essential for
Theorems \ref{lem1}, \ref{lem1+} and \ref{17.02.11} as
well as for Corollary \ref{co2}. Moreover, in this
example, the set $J$ defined in \eqref{setJ} consists
only of one point $2$.
   \begin{exa} \label{UniversalPictures}
Let $\varkappa=2$, $\vartheta_1 = 1$,
$\vartheta_2=\sqrt{a}$ and $\vartheta_3 = a$ with $a
\in (1,\infty)$. Then $\alpha = \frac{1}{\sqrt{a}}$,
$\alpha^\dag = \beta^\dag = 1$, $\beta =
\sqrtx{1}{4}{a}$ and $\gamma=\sqrt{a}$. Set
   \begin{align} \label{kopia}
   a_n = ((\alpha^\dag)^n + \gamma^n)^2 =
\vartheta_1^n + 2\vartheta_2^n + \vartheta_3^n, \quad
n \in \zbb_+.
   \end{align}
Clearly $\{a_n\}_{n=0}^\infty$ and
$\{\sqrt{a_n}\}_{n=0}^\infty$ are determinate
Stieltjes moment sequences. It follows from
\eqref{kopia} that $\supp \mu = \{\vartheta_1,
\vartheta_2, \vartheta_3\}$ and $\supp \nu =
\{\alpha^\dag, \gamma\}$. Hence $\frac{\gamma}{\beta}
\, \alpha < \alpha^\dag$, $\nu((\alpha, \beta)) > 0$
and $\alpha, \beta \not\in \supp\nu$. Moreover, we
have $\iota_{\mathrm{s}}=3 > 2 =
\iota_{\mathrm{s}}^*$.

We will show that $\{\sqrtk{a_n}\}_{n=0}^\infty$ is
not a Stieltjes moment sequence for every integer
$\varkappa \Ge 3$. Suppose that, contrary to our
claim, $\{\sqrtk{a_n}\}_{n=0}^\infty$ is a Stieltjes
moment sequence for some integer $\varkappa \Ge 3$.
Denote by $\nu_\varkappa$ the representing measure of
$\{\sqrtk{a_n}\}_{n=0}^\infty$. By \eqref{compsup} and
Corollary \ref{suppmu-cor}, $\card(\supp\nu_\varkappa)
< \infty$ and $\sqrtk{\vartheta_3} = \sup \supp
\nu_\varkappa$. Hence, in view of Theorem
\ref{suppmu}, $\supp \mu = \pi_\varkappa\big((\supp
\nu_\varkappa)^\varkappa\big)$. If $\supp\nu_\varkappa
= \{\sqrtk{\vartheta_3}\}$, then
$\supp\mu=\{\vartheta_3\}$, a contradiction.
Otherwise, by \eqref{porfelTer}, there exists $x\in
\supp\nu_\varkappa \cap (0,\sqrtk{\vartheta_3})$.
Then, by \eqref{porfelTer} again, $x^0
\big(\sqrtk{\vartheta_3}\, \big)^\varkappa, x^1
\big(\sqrtk{\vartheta_3}\, \big)^{\varkappa-1},
\ldots, x^\varkappa \big(\sqrtk{\vartheta_3}\,
\big)^0$ is a strictly decreasing sequence of
$\varkappa+1$ elements of $\supp\mu$. Since
$\varkappa+1 \Ge 4$, we arrive at a contradiction with
$\card(\supp\mu)=3$.
   \end{exa}
The subsequent example is related to Theorems
\ref{lem1+} and \ref{17.02.11}.
   \begin{exa} \label{UniversalPictures2}
Let $\varkappa=2$. Set $\vartheta_1 = \frac 12$,
$\vartheta_2=1$ and $\vartheta_3=9$. Then $\alpha =
\frac{1}{6}$, $\alpha^\dag = \frac 13$, $\beta^\dag =
\frac{1}{\sqrt{2}}$, $\beta = 1$ and $\gamma=3$. Set
$a_n = (\alpha^n + (\alpha^\dag)^n + \beta^n +
\gamma^n)^2$ for $n \in \zbb_+$. The Stieltjes moment
sequences $\{a_n\}_{n=0}^\infty$ and
$\{\sqrt{a_n}\}_{n=0}^\infty$ are determinate. It is
easily seen that $\supp\mu = \{\frac{1}{36},
\frac{1}{18}, \frac{1}{9}, \frac{1}{6}, \frac{1}{3},
\frac{1}{2}, 1, 3, 9\}$,
$\mu((\vartheta_1,\vartheta_2))=0$, $\{\vartheta_1,
\vartheta_2, \vartheta_3\} \subseteq \supp \mu$, $\nu
((\alpha,\beta)) > 0$ and $\{\alpha,\beta\} \subseteq
\supp \nu$. Moreover, we have
   \begin{align} \label{Jerzy}
\alpha^\dag < \beta^\dag, \quad \frac{\gamma}{\beta}
\, \alpha > \alpha^\dag, \quad \iota_{\mathrm{s}}=2 <
4= \iota_{\mathrm{s}}^*.
   \end{align}
   \end{exa}
Now we show that Theorems \ref{lem1+}\,(iii) and
\ref{17.02.11} do not cover all possible situations in
which supports of representing measures of the
$\varkappa\,$th roots of Stieltjes moment sequences
may be involved.
   \begin{exa} \label{UniversalPictures3}
Let $\varkappa=2$. Set $\vartheta_1 = 1$,
$\vartheta_2=a^2$ and $\vartheta_3=a^8$ with $a \in
(1,\infty)$. Then $\alpha = \frac{1}{a^4}$,
$\alpha^\dag = \frac {1}{a^2}$, $\beta^\dag = 1$,
$\beta = a$ and $\gamma=a^4$. Set $a_n = (\alpha^n +
\beta^n + \gamma^n)^2$ for $n \in \zbb_+$. Clearly
$\{a_n\}_{n=0}^\infty$ and
$\{\sqrt{a_n}\}_{n=0}^\infty$ are determinate
Stieltjes moment sequences. We verify directly that
$\supp\mu = \{\frac{1}{a^8}, \frac{1}{a^3}, 1, a^2,
a^5, a^8\}$, $\mu((\vartheta_1,\vartheta_2))=0$,
$\{\vartheta_1, \vartheta_2, \vartheta_3\} \subseteq
\supp \mu$, $\nu ((\alpha,\beta)) = 0$ and $\{\alpha,
\beta\} \subseteq \supp \nu$. Moreover, \eqref{Jerzy}
is satisfied.
   \end{exa}
We conclude this section with an example of a
Stieltjes moment sequence whose $\varkappa\,$th root
is a Stieltjes moment sequence for $\varkappa=2,4$,
but not for $\varkappa=3$ (consult Corollary
\ref{co2}).
   \begin{exa} \label{prznotsq}
Fix a real number $a > 1$. Set $\hat\alpha=a^{-33}$,
$\hat\beta=1$, $\hat\gamma=a^{3}$ and
   \begin{align*}
a_n = (\hat\alpha^n + \hat\beta^n + \hat\gamma^n)^4,
\quad n \in \zbb_+.
   \end{align*}
It is clear that $\{a_n\}_{n=0}^\infty$,
$\{\sqrt{a_n}\}_{n=0}^\infty$ and
$\{\sqrtx{2}{4}{a_n}\}_{n=0}^\infty$ are determinate
Stieltjes moment sequences whose representing measures
have finite supports. We claim that
$\{\sqrtx{2}{3}{a_n}\}_{n=0}^\infty$ is not a
Stieltjes moment sequence.

For this, denote by $\mu$ and $\nu_4$ the representing
measures of Stieltjes moment sequences
$\{a_n\}_{n=0}^\infty$ and
$\{\sqrtx{2}{4}{a_n}\}_{n=0}^\infty$, respectively. By
Theorem \ref{suppmu}, we have
   \begin{align*}
\supp \mu = \pi_4\big((\supp \nu_4)^4\big) =
\pi_4\big(\{\hat\alpha, \hat\beta, \hat\gamma\}^4\big)
= \{a^{3j-33i}\colon i,j \in \zbb_+, \, i+j \Le 4\}.
   \end{align*}
This implies that
   \begin{align} \label{15}
\card(\supp \mu) = 15.
   \end{align}
Set $\vartheta_1 = \hat\alpha \hat\gamma^{\varkappa -
1}$, $\vartheta_2=\hat\beta^\varkappa$ and
$\vartheta_3=\hat\gamma^\varkappa$ with $\varkappa=4$.
Clearly $\vartheta_1=a^{-24}$, $\vartheta_2=1$ and
$\vartheta_3=a^{12}$. Since $\supp\nu_4 =
\{\hat\alpha, \hat\beta, \hat\gamma\}$,
$\nu_4\big((\hat\alpha,\hat\beta)\big)=0$ and
$\hat\alpha \hat\gamma^3 < \hat\beta^4$, we infer from
Theorem \ref{lem1} that $\mu\big((\vartheta_1,
\vartheta_2)\big)=0$, $\vartheta_3=\sup \supp \mu$ and
$\vartheta_1, \vartheta_2, \vartheta_3 \in \supp\mu$.

Suppose that, contrary to our claim,
$\{\sqrtx{2}{3}{a_n}\}_{n=0}^\infty$ is a Stieltjes
moment sequence. Let $\nu_3$ be its representing
measure. In view of Corollary \ref{suppmu-cor},
$\supp\nu_3$ is finite. Hence, by Theorem
\ref{suppmu}, we have
   \begin{align}   \label{sup3}
\supp \mu = \pi_3\big((\supp \nu_3)^3\big).
   \end{align}
Let $\alpha$, $\beta$, $\gamma$, $\alpha^\dag$ and
$\beta^\dag$ be as in \eqref{theta2.5} and
\eqref{abplus} with $\varkappa=3$. Then
$\alpha=a^{-32}$, $\beta=1$, $\gamma=a^4$,
$\alpha^\dag=\beta^\dag=a^{-8}$ and
   \begin{align} \label{need}
\alpha\gamma < \alpha\gamma^2 = \vartheta_1 <1.
   \end{align}
By Theorem \ref{lem1+}\,(iii-a),
$\nu_3\big((\alpha,\beta)\big)=0$. It follows from
Corollary \ref{suppmu-cor}\,(iii) that $\gamma =
\sup\supp\nu_3$. Applying Theorem \ref{lem1}, we
deduce that $\alpha, \beta, \gamma \in \supp \nu_3$.
If $\card(\supp\nu_3)=3$, then $\supp\nu_3 = \{\alpha,
\beta,\gamma\}$, and by \eqref{sup3} we have
$\card(\supp \mu) = 10$, which contradicts \eqref{15}.
The remaining possibility is that
$\card(\supp\nu_3)\Ge 4$. Then there exists $x \in
\supp\nu_3 \cap [(0,\alpha) \cup (\beta, \gamma)]$. We
will show that $\card(\supp\mu) \Ge 16$, again
contradicting \eqref{15}.

Consider first the case of $x \in (0, \alpha)$. If
$\alpha = x \gamma$, then $x=a^{-36}$ and the
sequence\footnote{\;We use \eqref{sup3} to justify
that the terms of the sequences being considered are
in $\supp \mu$.}
   \begin{align*}
\{x^3, x^2\alpha, x\alpha^2, \alpha^3, x^2, x\alpha,
x\alpha\gamma, \alpha^2\gamma, x, x\gamma, x \gamma^2,
\alpha\gamma^2, 1, \gamma, \gamma^2, \gamma^3\}
\subseteq \supp \mu
   \end{align*}
is strictly increasing. If $\alpha < x \gamma$, then
\eqref{need} implies that the sequence
   \begin{align*}
\{x^3, x^2\alpha, x^2, x\alpha, x^2\gamma,
x\alpha\gamma, x, \alpha, x\gamma, \alpha \gamma,
x\gamma^2, \alpha\gamma^2, 1, \gamma, \gamma^2,
\gamma^3\} \subseteq \supp \mu
   \end{align*}
is strictly increasing as well. Finally, if $\alpha >
x \gamma$, then, by \eqref{need} again, the sequence
   \begin{align*}
\{\xi_n\}_{n=1}^{15}:=\{x^3, x^2\alpha, x^2,
x^2\gamma, x\alpha, x\alpha\gamma, x, x\gamma,
x\gamma^2, \alpha\gamma, \alpha\gamma^2, 1, \gamma,
\gamma^2, \gamma^3\} \subseteq \supp \mu
   \end{align*}
is strictly increasing and $\xi_8 < \alpha <
\xi_{10}$. If $\alpha \neq \xi_{9}$, then evidently
$\card(\supp\mu) \Ge 16$. Otherwise $x = a^{-40}$, and
thus $\xi_6 < \alpha^2 < \xi_7$, which yields
$\card(\supp\mu) \Ge 16$.

Let now $x \in (\beta,\gamma)$. Then by \eqref{need}
the sequence
   \begin{align*}
\{\alpha^3, \alpha^2, x\alpha^2, \alpha^2\gamma,
\alpha, x\alpha, x^2\alpha, x\alpha\gamma,
\alpha\gamma^2, 1, x, x^2, x^3, x^2\gamma, x\gamma^2,
\gamma^3\} \subseteq \supp \mu
   \end{align*}
is strictly increasing. This completes the proof of
our claim.
   \end{exa}
   \section{Square roots}
In this section, we concentrate on square roots of
Stieltjes moment sequences which have representing
measures supported in finite sets. For $M,N \in \nbb$,
we define the following classes of Stieltjes moment
sequences:
   \begin{enumerate}
   \item [$\bullet$] $\sfr_M$ stands for the set of
all Stieltjes moment sequences having representing
measures $\mu$ such that $\supp \mu \subseteq
(0,\infty)$ and $\card(\supp\mu) = M$,
   \item[$\bullet$]  $\ssqrt{M}$  stands for
the set of all sequences $\{a_n\}_{n=0}^\infty \in
\sfr_M$ such that $\{\sqrt{a_n}\,\}_{n=0}^\infty$ is a
Stieltjes moment sequence,
   \item[$\bullet$] $\ssqrtt{M}{N}$ stands for the
set of all $\{a_n\}_{n=0}^\infty \in \sfr_M$ such that
$\{\sqrt{a_n}\,\}_{n=0}^\infty \in \sfr_N$.
   \end{enumerate}
By \eqref{compsup}, any member of $\sfr_M$ is
determinate. It follows from Theorem \ref{suppmu} that
   \begin{align} \label{sum-a}
\ssqrt{M} = \bigcup_{N = 1}^\infty \ssqrtt{M}{N}.
   \end{align}
We now describe all pairs $(M,N)$ for which the
classes $\ssqrtt{M}{N}$ are nonempty.
   \begin{thm} \label{tw2}
Let $M,N \in \nbb$. Then the following conditions are
equivalent\/{\em :}
   \begin{enumerate}
   \item[(i)] $\ssqrtt{M}{N} \neq \varnothing$,
   \item[(ii)] $2 N -1 \Le M \Le \binom{N+1}{2}$,
   \item[(iii)] $\nfr_{-}(M) \Le N \Le \nfr_{+}(M)$,
   \end{enumerate}
where $\nfr_{-}(M) =
\big\lceil\frac{\sqrt{8M+1}-1}{2}\big\rceil$ and
$\nfr_{+}(M) = \big\lfloor\frac{M+1}{2}\big\rfloor$.
Moreover, the following holds\/{\em :}
   \begin{enumerate}
   \item[(a)] $\{\nfr_{-}(M)\}_{M=1}^\infty$
and $\{\nfr_{+}(M)\}_{M=1}^\infty$ are monotonically
increasing sequences,
   \item[(b)] for each $k \in \nbb$, exactly $k$
terms of $\{\nfr_{-}(M)\}_{M=1}^\infty$ are equal to
$k$,
   \item[(c)] for each $k \in \nbb$, exactly $2$
terms of $\{\nfr_{+}(M)\}_{M=1}^\infty$ are equal to
$k$.
   \end{enumerate}
   \end{thm}
   \begin{proof}
(i)$\Rightarrow$(ii) Take $\{a_n\}_{n=0}^\infty \in
\ssqrtt{M}{N}$. Let $\mu$ and $\nu$ be as in
\eqref{zal} with $\varkappa=2$. Then
$\supp\nu=\{\xi_1, \ldots, \xi_N\}$, where $0 < \xi_1
< \ldots < \xi_N$. First note that $2N -1 \Le M$.
Indeed, this can be inferred from \eqref{porfelTer}
and the following inequalities
   \begin{align*}
\xi_1^2 < \xi_1\xi_2 < \xi_2^2 < \ldots < \xi_{N-1}^2
< \xi_{N-1} \xi_N < \xi_N^2.
   \end{align*}
Set $\varOmega_N=\big\{(k,l)\in J_N \times J_N \colon
k\Le l\big\}$ with $J_N:=\{1, \ldots, N\}$. By
\eqref{porfelTer} the mapping $\varOmega_N \ni (i,j)
\longmapsto \xi_i\xi_j \in \supp\mu$ is surjective,
and thus $M \Le \card(\varOmega_N)=\binom{N+1}{2}$.

(ii)$\Rightarrow$(i) We shall prove that for every
$N\in \nbb$ and for every $M\in \nbb$ such that $2 N-1
\Le M \Le \binom{N+1}{2}$ there exists a sequence
$\xi_1 < \ldots < \xi_N$ of positive real numbers such
that\footnote{\;For simplicity, we write
$\big\{\xi_i\xi_j\colon i, j \in J_N\big\}$ in place
of $\big\{x \in \rbb | \; \exists i, j \in J_N \colon
x=\xi_i\xi_j\big\}$.}
$\card\big(\big\{\xi_i\xi_j\colon i, j \in
J_N\big\}\big) = M$. Once this is done, we see that
$\big\{(\sum_{j=1}^N \xi_j^n)^2 \big\}_{n=0}^\infty
\in \ssqrtt{M}{N}$.

We proceed by induction on $N$. The cases of $N=1$ and
$N=2$ are easily seen to hold. Suppose that our claim
is valid for a fixed integer $N$ which is greater than
or equal to $2$. Take $M \in \nbb$ such that
   \begin{align} \label{n+1ind}
2 (N+1)-1 \Le M \Le \binom{N+2}{2}.
   \end{align}

First we consider the case when
   \begin{align}   \label{m3n}
M \Le 3N.
   \end{align}
It is clear that $k_0:=\binom{N+2}{2} - 2(N+1) \Ge 0$
(because $N \Ge 2$). Set $k=M - 2(N+1)$. It follows
from \eqref{n+1ind} that $-1 \Le k \Le k_0$. From
\eqref{m3n} we infer that $k+3 \Le N+1$. Fix $\xi_1
\in (0,\infty)$ and $t \in (1,\infty)$, and set $\xi_j
= \xi_1 t^{k+j}$ for $j = 2, \ldots, N+1$. It is
easily seen that the sets
$\Xi_{N+1}:=\{\xi_i\xi_j\colon i,j \in J_{N+1}\}$ and
$\{\xi_1^2\} \cup \{\xi_1^2 t^{i}\}_{i=k+2}^{k+N+1}
\cup \{\xi_1^2 t^{i}\}_{i=2k+4}^{2k+2(N+1)}$ coincide.
Thus, since $k+3 \Le N+1$, we can arrange the elements
of the set $\Xi_{N+1}$ as follows
   \begin{align*}
\xi_1^2 < \xi_1 \xi_2 < \ldots < \xi_1 \xi_{k+3} <
\xi_2^2 < \xi_2\xi_3 < \xi_3^2 < \ldots < \xi_{N}^2 <
\xi_{N}\xi_{N+1} < \xi_{N+1}^2.
   \end{align*}
Hence $\card(\Xi_{N+1}) = M$.

Now we consider the remaining possibility, namely that
   \begin{align}   \label{m3n-}
M > 3N.
   \end{align}
Set $M^\prime = M - (N+1)$. Then by \eqref{n+1ind} we
have
   \begin{align*}
2N -1 = 3N - (N+1) \overset{\eqref{m3n-}}< M^\prime
\overset{\eqref{n+1ind}} \Le \binom{N+2}{2} - (N+1) =
\binom{N+1}{2}.
   \end{align*}
By the induction hypothesis applied to $M^\prime$,
there exists a sequence $\xi_2 < \ldots < \xi_{N+1}$
of positive real numbers such that
$\card\big(\big\{\xi_i \xi_j \colon i,j=2, \ldots,
N+1\big\}\big) = M^\prime$. Then there exists $\xi_1
\in (0,\xi_2)$ such that $\xi_1^2 < \xi_1 \xi_2 <
\ldots < \xi_1\xi_{N+1} < \xi_2^2$. Since $\xi_2^2 \Le
\xi_i\xi_j$ for all $i,j=2, \ldots, N+1$, we conclude
that $\card\big(\big\{\xi_i \xi_j \colon i,j=1,
\ldots, N+1\big\}\big) = M^\prime + N+1 = M$. This
completes the induction argument. Hence (i) is valid.

It is a matter of routine to show that the conditions
(ii) and (iii) are equivalent. The assertions (a) and
(c) are easily seen to hold, so we only explicitly
prove (b). Set $M_k=\binom{k+1}{2}$ for $k\in \nbb$.
Then $M_1=1$ and $M_{k+1} - M_k=k+1$ for $k\in \nbb$.
Note that $\nfr_{-}(1)=1$ and
   \begin{align*}
\frac{\sqrt{8M_k+1}-1}{2} = \frac{(2k+1)-1}{2}=k,
\quad k \in \nbb.
   \end{align*}
This implies that for every $k\in \nbb$ and for every
$M \in \nbb$ such that $M_k + 1 \Le M \Le M_{k+1}$,
$\nfr_{-}(M) = k+1$. This completes the proof.
   \end{proof}
Using assertions (a), (b) and (c) of Theorem
\ref{tw2}, one can easily specify successive terms of
the sequences $\{\nfr_{-}(M)\}_{M=1}^\infty$ and
$\{\nfr_{+}(M)\}_{M=1}^\infty$. Below, we list the
first fifteen terms of each of these sequences.
   \vspace{1ex}
   \begin{align*}
   \begin{tabular}{|r||c|c|c|c|c|c|c|c|c|c|c|c|c|c|c|c|c|}
\hline $M$ & $1$& ${\boldsymbol 2}$ & $3$ &
${\boldsymbol 4}$ & $5$ & $6$ & $7$ & $8$ & $9$ & $10$
& $11$ & $12$ & $13$ & $14$ & $15$ & $\cdots$
   \\ \hline  \hline
$\nfr_{-}(M)$ & $1$ & ${\boldsymbol 2}$ & $2$ &
${\boldsymbol 3}$ & $3$ & $3$ & $4$ & $4$ & $4$ & $4$
& $5$ & $5$ & $5$ & $5$ & $5$ & $\cdots$
   \\ \hline
$\nfr_{+}(M)$ & $1$ & ${\boldsymbol 1}$ & $2$ &
${\boldsymbol 2}$ & $3$ & $3$ & $4$ & $4$ & $5$ & $5$
& $6$ & $6$ & $7$ & $7$ & $8$ & $\cdots$
   \\ \hline
   \end{tabular}
   \end{align*}
   \begin{center} {\sf Table 1} \end{center}
   \vspace{1ex}
Applying \eqref{sum-a} and Theorem
\ref{tw2} (see also Table 1), we get the following
corollary.
   \begin{cor} \label{discr}
If $M\in\{2,4\}$, then $\ssqrt{M} = \varnothing$. If
$M \in \nbb \setminus \{2,4\}$, then the set $A_M:=\{N
\in \nbb\colon \nfr_{-}(M) \Le N \Le \nfr_{+}(M)\}$ is
nonempty, $\ssqrtt{M}{N} \neq \varnothing$ for every
$N \in A_M$ and $\ssqrtt{M}{N} = \varnothing$ for
every $N \in \nbb \setminus A_M$.
   \end{cor}
It is worth mentioning that the Stieltjes moment
sequences $\{a_n\}_{n=0}^\infty$ constructed in
Examples \ref{UniversalPictures},
\ref{UniversalPictures3} and \ref{UniversalPictures2}
belong to the classes $\ssqrtt{3}{2}$, $\ssqrtt{6}{3}$
and $\ssqrtt{9}{4}$, respectively. According to
Corollary \ref{discr}, the square root of a Stieltjes
moment sequence whose representing measure is
concentrated on either a two point or a four point
subset of $(0,\infty)$ is never a Stieltjes moment
sequence. Hence, the following holds.
   \begin{cor}[\mbox{\cite[Lemma 3.3]{jsjs}}]
If $\alpha_1$, $\alpha_2$,
$\vartheta_1$ and $\vartheta_2$ are positive real
numbers, then the sequence $\big\{\sqrt{\alpha_1
\vartheta_1^n + \alpha_2
\vartheta_2^n}\big\}_{n=0}^\infty$ is a Stieltjes
moment sequence if and only if
$\vartheta_1=\vartheta_2$.
   \end{cor}
Let $\{a_n\}_{n=0}^\infty$ be a Stieltjes moment
sequence whose representing measure is concentrated on
a three point subset of $(0,\infty)$. Then there exit
$\alpha_1, \alpha_2, \alpha_3, \vartheta_1,
\vartheta_2, \vartheta_3 \in (0,\infty)$ such that
$\vartheta_1 < \vartheta_2 < \vartheta_3$ and $a_n =
\alpha_1 \vartheta_1^n + \alpha_2 \vartheta_2^n +
\alpha_3 \vartheta_3^n$ for all $n \in \zbb_+$. Hence,
in view of \eqref{sum-a} and Corollary \ref{discr},
the sequence $\{\sqrt{a_n}\}_{n=0}^\infty$ is a
Stieltjes moment sequence if and only if
$\vartheta_2^2 = \vartheta_1 \vartheta_3$ and
$\alpha_2^2 = 4 \alpha_1 \alpha_3$.
   \section{Appendix}
Here we present a proof of the implication
(a)$\Rightarrow$(b) of Corollary \ref{cor1}\,(iv)
which is independent of Theorems \ref{lem1} and
\ref{lem1+}. First, we state an auxiliary result.
   \begin{lem} \label{stilim}
If $\{a_n\}_{n=0}^\infty$ is a Stieltjes moment
sequence with a representing measure $\mu$ supported
in $[0,1]$, then $\lim_{k \to \infty} a_k =
\mu(\{1\})$ and $\{a_n - \lim_{k \to \infty}
a_k\}_{n=0}^\infty$ is a Stieltjes moment sequence.
   \end{lem}
   \begin{proof}
By Lebesgue's monotone convergence theorem, $\lim_{k
\to \infty} a_k = \mu(\{1\})$, and thus $\{a_n -
\lim_{k \to \infty} a_k\}_{n=0}^\infty$ is a Stieltjes
moment sequence with a representing measure
$\mu|_{[0,1)}$.
   \end{proof}
   \begin{proof}[Proof of implication
{\rm (a)}$\Rightarrow${\rm (b)} of Corollary
\ref{cor1}\,{\rm (iv)}] According to \eqref{compsup}
and \eqref{esssup}, the Stieltjes moment sequences
$\{a_n\}_{n=0}^\infty$ and
$\{\sqrt[\leftroot{1}\uproot{4}\varkappa]{a_n}\,
\}_{n=0}^\infty$ are determinate. Without loss of
generality we can assume that $\vartheta_2=1$
(consider $\{\vartheta_2^{-n} a_n\}_{n=0}^\infty$
instead of $\{a_n\}_{n=0}^\infty$). We infer from
\eqref{esssup} applied to $\{a_n\}_{n=0}^\infty$ that
   \begin{align} \label{wz1}
\lim_{n \to \infty}
(\sqrt[\leftroot{1}\uproot{4}\varkappa]{a_n})^{1/n} =
1.
   \end{align}
By \eqref{esssup}, $\supp \nu \subseteq [0,1]$. It
follows from Lemma \ref{stilim} that $\lim_{n\to
\infty} \sqrt[\leftroot{1}\uproot{4}\varkappa]{a_n} =
\sqrt[\leftroot{1}\uproot{4}\varkappa]{\mu(\{1\})}$.
Now, by Lemma \ref{stilim} applied to
$\{\sqrt[\leftroot{1}\uproot{4}\varkappa]{a_n}\}_{n=0}^\infty$,
we see that
$\{\sqrt[\leftroot{1}\uproot{4}\varkappa]{a_n} -
\sqrt[\leftroot{1}\uproot{4}\varkappa]{\mu(\{1\})}\}_{n=0}^\infty$
is a Stieltjes moment sequence. According to our
assumption, $\{a_n - \mu(\{1\})\}_{n=0}^\infty$ is a
Stieltjes moment sequence with the representing
measure $\mu_1 = \mu|_{[0,\vartheta_1]}$. Since
$\vartheta_1 \in \supp \mu$, we get $\vartheta_1 =
\sup \supp \mu_1$, which together with \eqref{esssup}
leads to
   \begin{align} \label{wz2}
\lim_{n \to \infty} \big(a_n-\mu(\{1\})\big)^{1/n} =
\vartheta_1.
   \end{align}
By the mean value theorem, we have (recall that $a_n -
\mu(\{1\}) \Ge 0$ for all $n \Ge 0$)
   \begin{align*}
\sqrt[\leftroot{1}\uproot{4}\varkappa]{a_n} -
\sqrt[\leftroot{1}\uproot{4}\varkappa]{\mu(\{1\})}=
\frac{a_n - \mu(\{1\})} {\varkappa
\,\tau_n^{\frac{\varkappa - 1}\varkappa}}, \quad n =
0, 1, 2, \ldots,
   \end{align*}
where $\tau_n$ is a real number such that $\mu(\{1\})
\Le \tau_n \Le a_n$. This implies that
   \begin{align} \label{wz3}
\frac{a_n - \mu(\{1\})}{\varkappa
\,a_n^{\frac{\varkappa - 1}\varkappa}} \Le
\sqrt[\leftroot{1}\uproot{4}\varkappa]{a_n} -
\sqrt[\leftroot{1}\uproot{4}\varkappa]{\mu(\{1\})} \Le
\frac{a_n - \mu(\{1\})}{\varkappa
\mu(\{1\})^{\frac{\varkappa - 1}\varkappa}}, \quad n =
0, 1, 2, \ldots
   \end{align}
The conditions \eqref{wz1}, \eqref{wz2} and
\eqref{wz3} combined give
   \begin{align}   \label{wz4}
\lim_{n \to \infty}
\big(\sqrt[\leftroot{1}\uproot{4}\varkappa]{a_n} -
\sqrtk{\mu(\{1\})}\,\big)^{1/n} = \vartheta_1.
   \end{align}
This together with \eqref{esssup} implies that a
representing measure of the Stieltjes moment sequence
$\{\sqrt[\leftroot{1}\uproot{4}\varkappa]{a_n} -
\sqrt[\leftroot{1}\uproot{4}\varkappa]{\mu(\{1\})}\,\}_{n=0}^\infty$
is supported in $[0,\vartheta_1]$. Since
   \begin{align*}
\sqrt[\leftroot{1}\uproot{4}\varkappa]{a_n} =
\big(\sqrt[\leftroot{1}\uproot{4}\varkappa]{a_n} -
\sqrt[\leftroot{1}\uproot{4}\varkappa]{\mu(\{1\})}\,\big)
+ \sqrt[\leftroot{1}\uproot{4}\varkappa]{\mu(\{1\})},
\quad n=0,1,2, \ldots,
   \end{align*}
we deduce from \eqref{wz4} and $\mu(\{1\}) > 0$ that
the representing measure $\nu$ of
$\{\sqrt[\leftroot{1}\uproot{4}\varkappa]{a_n}\,\}_{n=0}^\infty$
is supported in $\big[0,\vartheta_1 \big] \cup \{1\}$
and $\{\vartheta_1,1\} \subseteq \supp \nu$.
   \end{proof}
   \bibliographystyle{amsalpha}
   
   \end{document}